%% file: main.tex
\documentclass[hidelinks,onefignum,onetabnum]{siamart220329}
\usepackage{lipsum}
\usepackage{amsfonts}
\usepackage{graphicx}
\usepackage{epstopdf}
\usepackage{mathrsfs}

\ifpdf
  \DeclareGraphicsExtensions{.eps,.pdf,.png,.jpg}
\else
  \DeclareGraphicsExtensions{.eps}
\fi

\newsiamremark{remark}{Remark}
\newsiamremark{hypothesis}{Hypothesis}
\crefname{hypothesis}{Hypothesis}{Hypotheses}
\newsiamthm{claim}{Claim}

\usepackage{array}
\usepackage{stix}
\usepackage{amsmath}

\usepackage{caption}
\usepackage{subcaption}

\headers{Polytope Division Method}{E. Nielen, O.T.C. Tse, K. Veroy}

\author{Evie Nielen \thanks{University of Technology, Eindhoven
  (\email{e.t.e.nielen@tue.nl}}
\and Oliver Tse \thanks{University of Technology, Eindhoven
  (\email{o.t.c.tse@tue.nl}}
\and Karen Veroy \thanks{University of Technology, Eindhoven 
  (\email{k.veroy@tue.nl}}}

\usepackage{amsopn}
\usepackage{bbm}
\usepackage{import}
\usepackage{xcolor}
\usepackage{algcompatible}
\usepackage{algorithm}
\usepackage[noend]{algpseudocode}
\floatname{algorithm}{Algorithm}

\usepackage{etoolbox}\AtBeginEnvironment{algorithmic}{\small} 
\usepackage{float}
\newfloat{algorithm}{H}{lop}
\usepackage{comment}
\usepackage{tikz}
\usepackage{subcaption}
\usepackage{hyperref}
\usepackage{enumitem}

\makeatletter
\def\BState{\State\hskip-\ALG@thistlm}
\makeatother

\newcommand{\eps}{\varepsilon}

\newcommand{\R}{\mathbb{R}}
\newcommand{\N}{\mathbb{N}}
\newcommand{\V}{\mathbb{V}}
\newcommand{\simp}{\mathbb{S}}
\newcommand{\bp}{\mathbb{B}_{P}}
\newcommand{\bI}{\mathbf{I}}
\newcommand{\bX}{\mathbf{X}}
\newcommand{\calB}{\mathscr{B}}
\newcommand{\calL}{\mathcal{L}}
\newcommand{\calP}{\mathcal{P}}

\newcommand{\calD}{\mathscr{D}}
\newcommand{\calS}{\mathscr{S}}
\newcommand{\calH}{\mathcal{H}}

\newcommand{\calN}{\mathcal{N}}

\newcommand{\calV}{\mathcal{V}}
\newcommand{\err}{J}

\newcommand{\RB}{\textup{RB}}

\newcommand{\conv}{\mathrm{Conv}}
\newcommand{\samsize}{|S|}

\newcommand{\FL}{\mathsf{FL}}
\newcommand{\bgamma}{\boldsymbol{\gamma}}
\newcommand{\etab}{\boldsymbol{\eta}}
\newcommand\norm[1]{\left\|#1\right\|}
\newcommand{\norml}[1]{{\left\vert\kern-0.25ex\left\vert\kern-0.25ex\left\vert #1\right\vert\kern-0.25ex\right\vert\kern-0.25ex\right\vert}_{\lambda}}
\newcommand{\interior}[1]{%
  {\text{int}(#1)}%
}
\newcommand{\sfb}{\mathsf{b}}
\newcommand{\sfP}{\mathsf{P}}

\DeclareMathOperator*{\argmin}{arg\,min}
\DeclareMathOperator*{\argmax}{arg\,max}

\title{Polytope Division Method: A Scalable Sampling Method for problems with high-dimensional parameters \thanks{\funding{This work is part of a project that has received funding from the European Research Council (ERC) under the
European Union’s Horizon 2020 Research and Innovation Programme (Grant Agreement No. 818473).}}}
\date{\today}

\begin{document}

\maketitle
\begin{abstract}
Configuration Optimization Problems (COPs), which involve minimizing a loss function over a set of discrete points $\bgamma\subset P$, are common in areas like Model Order Reduction, Active Learning, and Optimal Experimental Design. While exact solutions are often infeasible, heuristic methods such as the Greedy Sampling Method (GSM) provide practical alternatives, particularly for low-dimensional cases. GSM recursively updates 
$\bgamma$ by solving a continuous optimization problem, which is typically approximated by a search over a discrete sample set $S\subset P$. However, as the dimensionality grows, the sample size suffers from the curse of dimensionality.

To address this, we introduce the Polytope Division Method (PDM), a scalable greedy-type approach that adaptively partitions the parameter space and targets regions of high loss. PDM achieves linear scaling with problem dimensionality and offers an efficient solution approach for high-dimensional COPs, overcoming the limitations of traditional methods.
\end{abstract}



\tableofcontents

\section{Introduction}
\label{sec: Introduction}
\input{Introduction}

\section{Greedy Sampling Method}
\label{sec: GSM}
\input{GreedySamplingMethod}

\section{Polytope Division Method}
\label{sec: PDM}
\input{PolytopeDivisionMethod}

\section{COPs in Model Order Reduction}
\label{sec: cops in mor}
In this Section, we look at different examples in Model Order Reduction and formulate it as a COP.
\subsection{Reduced Basis Method}
\label{subsec: Reduced Basis Method}
\input{Reduced_Basis_Method}

\subsection{Empirical Interpolation Method}
\label{subsec: Empirical Interpolation Method}
\input{Empirical_Interpolation_Method}

\section{Numerical Tests}
\label{sec: Numerical examples}
\input{Numerical_examples}

\section{Outlook}
\label{sec: Outlook}
\input{Outlook}

\appendix
\section{Proof of Lemma \ref{lemma: facets p-boundary polytopes}}
\label{app: proofs}
\input{Proofs}

\section*{Acknowlegdements}
The authors thank Pierre Mollo and Rahul Dhopeshwar for their help with the PDM code. The authors acknowledge using Grammarly and chatGPT to polish the written text for spelling, grammar, and general style.
\bibliographystyle{siamplain}

\input{main.bbl}
\end{document}

%% file: Introduction.tex
Optimization problems are ubiquitous in computational science, engineering, and machine learning. In this work, we focus on a specific class known as \emph{Configuration Optimization Problems} (COPs), where the objective is to minimize a loss function that depends on configuration of points. Such problems can generally be stated as follows: For a fixed number of points $n\in\N$ in a compact set, $P\subset\R^d$, and a given \emph{loss function}, $\calL\colon \varGamma(P)\to [0,+\infty)$, we want to find
\begin{align}\label{problem: global-optimization}
    \bgamma \in \argmin_{\etab\in \varGamma_n} \calL(\etab).\tag{COP}
\end{align}
Here, $\Gamma(P)$ is the space of finite configurations, defined by
\[
    \Gamma(P) \coloneq \bigsqcup_{n\in \N} \Gamma_n(P),
\]
where $\Gamma_n(P) \coloneq \{ \bgamma\subset P : \#\bgamma = n\}$ is the family of finite subsets of $P$ with cardinality $n$.
In other words, the task at hand is to find a subset $\bgamma \subset P$ of $n \in\N$ points that minimizes the loss function $\calL$.

It is, typically, computationally infeasible to solve these problems exactly, but strategies exist to heuristically find (possibly suboptimal) solutions to \eqref{problem: global-optimization}. One family of such strategies is the greedy method \cite{cohen2020reduced,hesthaven2014efficient,jiang2020adaptive,jiang2017offline,urban2014greedy,veroy2003posteriori,wu2019active}, which iteratively builds $\bgamma \subset P$ by choosing the (locally) optimal update at each iteration.
In this approach, $\bgamma^1 = \{p_1\}$ is chosen arbitrarily; then, in the $j$-th step, we update the configuration with $p_j$, i.e., $\bgamma^j=\bgamma^{j-1}\cup\{p_j\}$, where
\begin{align}
\label{eq: strong-greedy}
p_j \in \argmin_{q \in P\setminus \bgamma^{j-1}} \calL(\bgamma^{j-1} \cup \{q\}),\qquad \bgamma^{j-1}\in\varGamma_{j-1}.
\end{align}
While this procedure does not generally lead to a solution of \eqref{problem: global-optimization}, it does provide a tractable strategy to approximate solutions to \eqref{problem: global-optimization} by choosing the best possible update in the $j$-th step, given the configuration in the $(j\!-\!1)$-th step.

Unfortunately, \eqref{eq: strong-greedy} remains computationally intensive in many practical applications due to the existence of multiple local minimizers. To empirically overcome this, the Greedy Sampling Method (GSM) is often used. GSM aims to approximate the global minimizer of~\eqref{eq: strong-greedy} by determining a minimizer inside a sample set $S\subset P$. However, GSM suffers from two main problems:
(1) the lack of a sufficient understanding of the selection of method parameters, in particular, the sample size $\samsize$, and (2) scalability issues in higher dimensions or the "curse-of-dimensionality" (see Section~\ref{sec: GSM} for details).

In this paper, we propose an alternative greedy approach called the \emph{Polytope Division Method} (PDM) that is capable of circumventing the aforementioned pitfalls of GSM. The key idea behind PDM---based on a divide-and-conquer strategy---is to recursively divide the parameter set into polytopes where, in each step, polytope configurations are further divided in accordance with the updated configuration $\bgamma$.
The number of samples required to find an approximate solution $\bgamma \in \varGamma_n$ for fixed $n\ge 1$ scales only linearly with the dimension. Additionally, the number of samples is updated at each step, eliminating the need to choose a sample size $\samsize$ beforehand.

\subsection*{COPs in Applications}
We now consider two examples of COPs.

\subsubsection*{Example 1: Reduced Basis Methods (RBM)}

The first example is taken from Reduced Basis Methods (RBM) in Model Order Reduction (MOR) for parametrized families of partial differential equations (PDEs). RBM is used in many applications, such as material engineering \cite{huynh2007reduced,ngoc2005certified}, electro-magnetism \cite{chen2010certified}, heat flow \cite{rozza2009reduced} and many more \cite{deparis2008reduced,negri2013reduced,rozza2011reduced,rozza2007stability,sen2006natural}. It is important to note that this is only a small selection of the extensive literature available on the topic. In RBM, we seek a set of basis vectors that provide good approximations to the solution manifold of parametrized PDEs. 

Consider a parameter set $P$ and a PDE with solution $u(p)\in \calV$ depending on the parameter $p \in P$. Then a \emph{reduced basis} made up of a family of solutions at different parameter points $p\in P$ is used to approximate the full solution manifold $\mathcal{M} \coloneq \{u(p)\,|\, p \in P\}$. In this context, the loss function can be written as 
\begin{align}
\label{eq: loss-reduced-basis}
\calL_\RB(\bgamma) = \max_{q \in P} \norm{u(q) - \sfP_{V_{\bgamma}} u(q)}^2_{\calV}.
\end{align}
Here, $V_{\bgamma} = \text{span}\{u(p) \,| \, p \in \bgamma\}$ is the linear span of solutions associated to a configuration of parameter points $\bgamma$, and $\sfP_{V_{\bgamma}}$ is the projection operator onto the reduced basis space $V_{\bgamma}$. More details are provided in Section~\ref{subsec: Reduced Basis Method}.

\subsubsection*{Example 2: Empirical Interpolation Method (EIM)} The second example is given by the Empirical Interpolation Method (EIM) as introduced in \cite{barrault2004} and further discussed, among others, in \cite{chen2014weighted}, \cite{eftang2010posteriori}, \cite{maday2013generalized}. EIM is used to find an approximate affine decomposition of a function $s\colon P\to \R^N$ for some $N \in \N$. In this case, the loss function reads
\begin{align}
\label{eq: loss-EIM}
\calL_{\text{EIM}}(\bgamma) = \max_{q \in P} \norm{s(q) - s_{\bgamma, \mathbf{I}}(q)}.
\end{align}
Here, $s_{\bgamma, \bI}$ is the EIM approximation, which depends on the set of configuration points $\bgamma$ and a set of interpolation points $\bI$. Details about the function $s_{\bgamma, \bI}$ are given in Section~\ref{subsec: Empirical Interpolation Method}.

\medskip

In this work, we mainly focus on the previous two examples. Still, more examples of this class of optimization problems can be found, for instance, in the context of Optimal Experimental Design (OED) \cite{ucinski2004optimal} and in active learning for regression \cite{wu2019active}, where our PDM approach applies. 

\subsection*{Outline of paper}
In Section~\ref{sec: GSM}, we provide a more detailed explanation of GSM, as well as a literary overview. In Section~\ref{sec: PDM}, we describe the Polytope Division Method in detail. In Section~\ref{sec: cops in mor}, we expand on the first two examples of COPs, the Reduced Basis Method, and the Empirical Interpolation Method. We test PDM in Section~\ref{sec: Numerical examples} with a thermal block example and a heat conductivity problem with a Gaussian heat source. Lastly, Section~\ref{sec: Outlook} presents an outlook for further research.

%% file: GreedySamplingMethod.tex
We begin this section by providing a more detailed summary of the Greedy Sampling method (GSM) \cite{veroy2003posteriori}.
The key idea behind GSM is to compute a heuristic solution to~\eqref{problem: global-optimization} by sequentially selecting points to construct the configuration $\bgamma$. Ideally, the configuration $\bgamma \in \varGamma_{j-1}$ should be updated with the point $p_j$ given by~\eqref{eq: strong-greedy}. However, even the task of determining this optimal update can be infeasible in many applications. Therefore, the optimizer of the \textit{continuous} optimization problem~\eqref{eq: strong-greedy} is approximated by the optimizer of a \textit{discrete} optimization problem over a sample set.

The precise implementation of GSM depends on the form of the loss function. We focus on the class of loss functions involving a maximum, i.e., where there exists an objective function such that $\calL(\bgamma) = \max_{p \in P} \err(p, \bgamma)$, where $\err(p, \bgamma)$ has the following properties
\begin{enumerate}
    \item $\err(p, \bgamma) \ge 0$,
    \item The function $\err(\cdot, \bgamma)$ is Lipschitz continuous with Lipschitz constant $L_\err$.
    \item For every $p \in \bgamma$, it holds that
    $$\err(p, \bgamma) = 0.$$
    \item For every $p \in P$, it holds that if $\bgamma^n \subset \bgamma^m$, then
    \begin{align*}
        \err(p, \bgamma^m) \le \err(p, \bgamma^n).
    \end{align*}
\end{enumerate}
In GSM, we then consider a finite sample set $S \subset P$ and, in each step, determine
\begin{align}
\label{eq: error-estimates}
p_j \in \argmax_{p \in S\setminus \bgamma^{j-1}} \err(p, \bgamma^{j-1}),
\end{align}
GSM is described by Algorithm~\ref{gs} (cf. \cite{quarteroni2015}).

\begin{algorithm}
\caption{Greedy Sampling Method}\label{gs}
\begin{algorithmic}[1]
\State Initialize tolerance {\tt tol}, choose sample size $\samsize$
\State Choose sample set $S \subset P$ of cardinality $\samsize$ 
\State {\tt err} $\gets$ {\tt tol} $+ 1$
\State Pick $p \in P$
\State Set $\bgamma \coloneq \{p\}$
\While{{\tt err}$ \, >$ {\tt tol}}
\State Determine $p \in \displaystyle\argmax_{q \in S\setminus\bgamma} \err(q, \bgamma)$
\State {\tt err} $\gets \displaystyle \err(p, \bgamma)$
\State $\bgamma \gets \bgamma \cup \{p\}$ 
\EndWhile
\end{algorithmic}
\end{algorithm}

Different variants of GSM exist based on the selection of the sample set $S$. This set could be constructed via random sampling inside the parameter set $P$; we refer to this variant as GSM-Random. Another option is to construct the set $S$ via Latin Hypercube Sampling (GSM-LHS). Latin Hypercube Sampling seeks to ensure a better representation of the variability in the parameter set; see  e.g. \cite{mckay2000comparison} for details.

GSM is a strategy to compute a heuristic, possibly suboptimal, solution to \ref{problem: global-optimization}, but the method has drawbacks:
\begin{enumerate}
    \item In GSM, the sample size $\samsize$ must be selected a priori, without clear guidelines for a suitable size.
    The sample size $\samsize$ should be chosen to ensure that the sample set $S$ adequately represents the set $P$. Nonetheless, we are unaware of any theoretical results that provide a priori guidance regarding the required sample size to guarantee a good approximation of the optimizer.
    \item Commonly used heuristic strategies choose $\samsize$ by taking, for example, $\samsize \approx n^d$, thus taking roughly $n$ parameter samples per parameter direction. However, such heuristic strategies rapidly face the "curse of dimensionality". Since the objective functions~\eqref{eq: error-estimates} have to be computed for each sample in every greedy step, the computations quickly become too costly in higher dimensions.
\end{enumerate}

The first drawback of GSM is addressed in \cite{cohen2020reduced}, specifically for the reduced basis setting. Here, sample sizes are given to guarantee convergence with high probability under the assumption that the solution $u(p)$ is analytic. However, this assumption is not met in general.

Several adaptations of GSM have been proposed to overcome the second drawback. Particularly in RB methods, the computational burden is decreased in \cite{ballani2016reduced} by formulating the problem in a low-rank tensor format resulting in more efficient computations. Other approaches to reduce the computational costs focus on the reduction of the sample sizes $\samsize$. In \cite{sen2008reduced}, GSM is performed sequentially over smaller, disjoint sample sets. In \cite{hesthaven2014efficient}, the sample set is updated after each step by discarding irrelevant samples and enriching the set with more relevant samples based on the value of the objective function. In \cite{jiang2017offline}, the authors consider a surrogate training set constructed via the successive maximization method. Hybrid versions of the aforementioned methods are investigated in \cite{jiang2020adaptive}. In \cite{eftang2012parameter}, alternative methods are proposed to reduce the high computational cost associated with large configurations in the context of the empirical interpolation method (EIM) \cite{barrault2004}. In particular, \cite{eftang2012parameter} proposed an $hp$-refinement procedure that subdivides the parameter domain and subsequently constructs smaller configurations for each subdomain.  In the EIM context, the smaller configurations result in a reduction in online cost, but comes at the price of a significantly more expensive offline stage \cite{eftang2012parameter}.  Furthermore, in each subdomain, the standard greedy sampling method (with its drawbacks) is still employed in each subdomain. Finally, \cite{urban2014greedy} uses gradient-based methods to compute local optima from a coarse set of starting points.
The maximizer of the largest local maximum is then selected to update the configuration $\bgamma$.

Although these alternatives lead to faster computations, they eventually suffer from the same scalability issues as the original GSM. To overcome this issue, we propose the Polytope Division Method (PDM), a greedy method that divides the parameter set into polytopes and uses the barycenters as samples. We introduce this method next.

%% file: PolytopeDivisionMethod.tex
\begin{figure}[h]
\centering
    \begin{subfigure}[t]{0.17\linewidth}
        \centering
        \input{Figures/S1_PDM-2d}
        \caption*{Step (1)}
        \label{fig: sf1: PDM2d}
    \end{subfigure}
    \hspace{.2cm}
    \begin{subfigure}[t]{0.17\linewidth}
        \centering
        \input{Figures/S2_PDM-2d}
        \caption*{Step (2)}
        \label{fig: sf2: PDM2d}
    \end{subfigure}
    \hspace{.2cm}
    \begin{subfigure}[t]{0.17\linewidth}
        \centering
        \input{Figures/S3_PDM-2d}
        \caption*{Step (3)}
        \label{fig: sf3: PDM2d}
    \end{subfigure}
    \hspace{.2cm}
    \begin{subfigure}[t]{0.17\linewidth}
        \centering
        \input{Figures/S4_PDM-2d}
        \caption*{Step (4)}
        \label{fig: sf4: PDM2d}
    \end{subfigure}
    \hspace{.2cm}
    \begin{subfigure}[t]{0.17\linewidth}
        \centering
        \input{Figures/S5_PDM-2d}
        \caption*{Step (5)}
        \label{fig: sf5: PDM2d}
    \end{subfigure}
\caption{Depiction of the steps in PDM for the $2$-dimensional parameter case. (1) Sample first parameter. (2) Compute barycenters. (3) Select the barycenter with the largest objective function value. (4) Update polytope division. (5) Compute the new barycenters.}
\label{fig: PDM2d}
\end{figure}
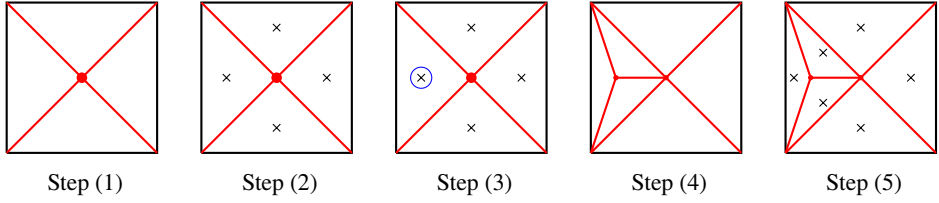

In this section, we introduce the Polytope Division Method. In a nutshell, the key idea behind the Polytope Division Method is to divide the parameter domain $P$ into polytopes and recursively compare the objective function values on the barycenters of these polytopes. After selecting the barycenter with maximal value and updating the configuration accordingly, we then refine the polytope division.
Figure~\ref{fig: PDM2d} depicts an example of the first steps of PDM in the $2$-dimensional case. 

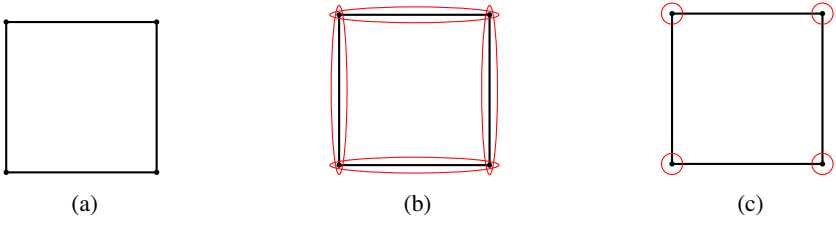
\begin{figure}[H]
\centering
    \begin{subfigure}[b]{0.32\linewidth}
        \centering
        \input{Figures/P}
        \caption{}
        \label{fig: P1}
    \end{subfigure}
    \hfill
    \begin{subfigure}[b]{0.32\linewidth}
        \centering
        \input{Figures/P_facets}
        \caption{}
        \label{fig: P2-facets}
    \end{subfigure}
    \hfill
    \begin{subfigure}[b]{0.32\linewidth}
        \centering
        \input{Figures/P_vertices}
        \caption{}
        \label{fig: P3-vertices}
    \end{subfigure}
\caption{Example of a $2$-dimensional polytope along with its facets and vertices. (a) Example of a polytope, (b) facets of a polytope, and (c) vertices of a polytope.}
\label{fig: pol_fac_ver}
\end{figure}

Before we explain the algorithm in detail, we first define the necessary terms: polytopes, facets and barycenters. The definition of the polytope requires the definition of the \emph{convex hull}
\cite{brondsted1983}.
\begin{definition}
The \textbf{convex hull} of any set $M \subset \R^d$ is the smallest convex set containing $M$. We denote the convex hull of $M$ by $\conv(M)$.
\end{definition}
\begin{definition}
    The \textbf{affine hull} of any $M \subset \R^d$, denoted by $\textup{aff}(M)$ is the smallest affine subspace containing $M$. For any affine subspace $A \subset \R^d$, there exists an $x \in \R^d$ and linear space $L$ such that $A = x + L$. The \textbf{dimension} of $A$, denoted by $\textup{dim}(A)$ equals the dimension of the corresponding linear space $L$.
\end{definition}
\begin{definition}[Polytopes, faces, facets, and vertices]
\label{def: polytope}
\begin{enumerate}[label={(\roman*)}]
    \item A \textbf{polytope} $P$ in $\R^d$ is the convex hull of a non-empty finite subset $M \subset \R^d$ of the form $M=\{x_1, \ldots, x_m\}$ with $m \ge d+1$, i.e., $P=\conv(M)$. The dimension of a polytope $P$ is given by the dimension of its affine hull, i.e., $\mathrm{dim}(\mathrm{aff}(P))$.

    \item A \textbf{face} $F$ of a polytope $P$ is a convex subset of $P$ such that there exists $c \in \R^d$ and $c_0$ satisfying
    \begin{itemize}
        \item[--] $F = P \cap \{x \in \R^d: c \cdot x = c_0\}$, and
        \item[--] for all $x \in P$ it holds that $c \cdot x \le c_0$.
    \end{itemize}

    \item A \textbf{facet} of a polytope $P$ is a face $F$ of $P$ with $0 \le \mathrm{dim}(F) = \mathrm{dim}(P) - 1$.
    \item A \textbf{vertex} $v$ of $P$ is a face $F$ of $P$ with $\mathrm{dim}(F) = 0$.
\end{enumerate}
\end{definition}
The definitions of faces, facets, and vertices are based on \cite{ziegler1994lectures}. An example of a polytope, its facets and its vertices is given in Figure~\ref{fig: pol_fac_ver}. Facets play an important role in PDM. We thus introduce a notation referring to facets.
\begin{definition}
    Let $P \subset \R^d$ be a polytope with facets $F_1, \ldots, F_M$. We define the set of facets of $P$ as $\partial P \coloneq \{F_i \;|\; i = 1, \ldots, M\}$.

    Moreover, $\partial_2 P$ denotes the set of facets of the facets of $P$. We refer to the elements in this set as the $2$-subfacets of $P$.
    Similarly, $\partial_k P$ contains the $k$-subfacets of $P$.
\end{definition}
A particular type of polytope is the simplex, which we define as
\begin{definition}
    $\sigma$ is a \textbf{$d$-simplex} if it is a $d$-dimensional polytope with precisely $d+1$ vertices.
    We denote the family of all simplices by $\simp$.
\end{definition}
Another particular type of polytope is the $P$-boundary polytope.
\begin{definition}
    Let $P$ be a $d$-dimensional polytope. $\tau$ is a \textbf{$P$-boundary polytope} if $\tau$ is a $d$-dimensional polytope and there exists an $n \in \N, n \le d-1$, points $x_1, \ldots, x_n \in P$ and some $F \in \partial_n P$ such that $\tau = \conv(\{x_1, \ldots, x_n\} \cup F)$.

    We denote the family of all $P$-boundary polytopes by $\bp$
\end{definition}
We remark that we do not include the case $n = d$ in the definition above. This allows us to distinguish between simplices and $P$-boundary polytopes. Indeed, notice that for $n = d$, the object $\conv(\{x_1, \ldots, x_n\} \cup F)$ with $F \in \partial_n P$ is a simplex.

Note that Definition~\ref{def: polytope} states that a polytope is a convex hull of points $\{x_1, \ldots, x_m\}$. In particular, a polytope $P$ can be expressed as the convex hull of its vertices $\{v_1, \ldots, v_n\}$, with $d+1\le n\le m$. With this representation, we can define the barycenter of a polytope.
\begin{definition}
    Let $P \subset \R^d$ be a polytope with vertices $v_1, \ldots, v_n$, in other words, $P = \conv(\{v_1, \ldots, v_n\})$. The \textbf{barycenter} of $P$ is then given by
    \[
    	\sfb_P \coloneq \frac{1}{n}\sum_{i=1}^n v_i. 
    \]
\end{definition}
In PDM, we assume that the parameter set $P$ is a hyperrectangle and we divide this into polytopes in the following sense.
\begin{definition}
    Let $P \subset \R^d$ be a $d$-dimensional polytope. Then a family of polytopes $\calD$ is a \textbf{proper polytope division} of $P$ if 
 \begin{enumerate}[label=(\roman*)]
 	\item $D \subset \R^d$, $D\in\calD$ are $d$-dimensional polytopes and $P = \bigcup_{D\in \calD} D$;
   \item For any $D,D'\in\calD$, $D\ne D'$, $\mathrm{int}(D)\cap \mathrm{int}(D') = \emptyset$.
 \end{enumerate}
\end{definition}

We can now use the proper polytope division and the barycenters to determine the update procedure of a configuration. Suppose that, at the $(j\!-\!1)$-th step of the configuration construction, we have a proper polytope division $\calD^{j-1}$ of $P$.
We then select
\begin{align}
\label{eq: D selection}
    D \in \underset{E\in\calD^{j-1}}{\arg\max}\, \,\err(\sfb_E, \bgamma^{j-1}) \tag{\#},
\end{align}
and set $\bgamma^j=\bgamma^{j-1}\cup \{\sfb_D\}$.
Recall that $\sfb_D$ denotes the barycenter of polytope $D$. \eqref{eq: D selection} thus chooses that polytope $D \in \calD$ that maximizes the objective function $\err$ (see Step 3 of Fig.~\ref{fig: PDM2d}), and we then proceed to add $\sfb_D$ to the configuration $\bgamma$ (see Step 4 of Fig.~\ref{fig: PDM2d}).

Now that we know how to update subspace $\bgamma^j$ for a given proper polytope division, the question arises of how to divide the hyperrectangle $P$ in the first place. A naïve approach would be to slice the parameter space at equidistant intervals in every direction. This leads to division into multiple smaller hyperrectangles. However, in this approach, the number of polytopes scales exponentially with respect to the dimension $d$, making it computationally unfeasible to evaluate the objective function at each barycenter. Thus, an alternative division is necessary.

\begin{figure}[h]
    \begin{subfigure}[b]{0.49\linewidth}
        \centering
        \input{Figures/FS1point-facet}
        \caption*{(a) Example of a point and a facet}
        \label{fig: FS1}
    \end{subfigure}
    \hfill
    \begin{subfigure}[b]{0.49\linewidth}
        \centering
        \input{Figures/FS2convex-hull}
        \caption*{(b) Convex hull of point and facet}
        \label{fig: FS2}
    \end{subfigure}
\caption{Example of the convex hull of a point and a facet. The dashed line in (a) represents the polytope, the red dot a point $p$, and the thick red line the facet $F$. The filled area in (b) represents $\conv(p \cup F)$.}
\label{fig: conv_point_facet}
\end{figure}
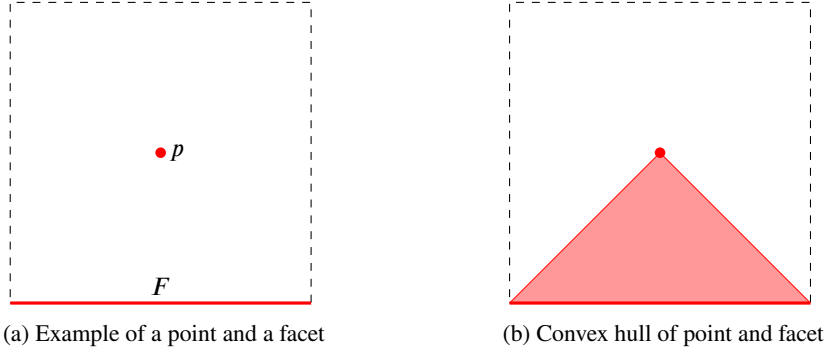

For a $d$-dimensional hyperrectangle, the number of facets equals $2d$, and this number, therefore, scales linearly with respect to the dimension. Consequently, it seems reasonable to propose a division using the facets of the original hyperrectangle. This is precisely the key idea behind PDM. In PDM, we use an operation called \emph{facet linking}, in which we link a point $p\in P$ to all the facets of the polytope $P$. Before we define facet linking, we explain the notation~$\conv(p \cup F)$, where $p \in F$ is a point and $F$ is a facet of $P$. The convex hull $\conv(p\cup F)$ is the smallest convex set that contains both $p$ and $F$. Figure \ref{fig: conv_point_facet} displays a $2$-dimensional example of such a convex hull.

\begin{definition}
    Let $P \subset \R^d$ be a polytope and $p \in P$  be an arbitrary point.
    The \textbf{facet linking operator} $\FL$ is given by
    \begin{align*}
        \FL(p, P) = \big\{\conv(p\cup F) \, \big|\, F\in\partial P\big\}.
    \end{align*}
\end{definition}
\noindent Note that $\FL(p, P)$ is a proper polytope division of $P$ for any $p \in \interior{P}$.

\medskip

Suppose that, in the $(j\!-1\!)$-th step of the algorithm, we have the configuration $\bgamma^{j-1} \in \varGamma_{j-1}$ and the proper polytope division $\calD^{j-1}$. Furthermore, suppose we update $\bgamma$ with $\sfb_D$ for some $D\in \calD^{j-1}$. The new polytope division of $P$  then becomes
\begin{align}
\calD^{j} = \bigl(\calD^{j-1} \backslash \{D\}\bigr) \cup \FL(\sfb_D, D).
\end{align}
The Polytope Division Method is summarized in Algorithm \ref{pdm}.

\begin{algorithm}
\caption{Polytope Division Method}\label{pdm}
\begin{algorithmic}[1]
\State Initialize tolerance {\tt tol}
\State {\tt err} $\gets$ {\tt tol} $+ 1$
\State Choose $p \in \interior{P}$
\State Set $\bgamma \coloneq \{p\}$
\State Set $\calD\coloneq\{D\in \FL(p, P)\}$
\While{{\tt err}$ \, >$ {\tt tol}}
\State Compute $D \in \displaystyle \arg\max_{E\in\calD}\err(\sfb_E, \bgamma)$
\State {\tt err} $\gets \err(\sfb_{D}, \bgamma)$
\State $\bgamma \gets \bgamma \cup \{\sfb_D\}$
\State $\calD \gets (\calD \backslash \{D\}) \cup \FL(\sfb_D,D)$
\EndWhile
\end{algorithmic}
\end{algorithm}

Figure \ref{fig: PDM2d} illustrates the steps of PDM only for the case $d=2$, but PDM works for general dimensions. Nonetheless, we must note the following: 
\begin{enumerate}
    \item \label{issue: shapes} For $d \ge 3$, the polytope division consists of more shapes than just simplices.
    \item \label{issue: facets} It is not trivial to determine the facets of a general polytope described by its vertices.
\end{enumerate}
Expanding on the first point, we note that, in the $2$-dimensional case, we only divide the polytopes into triangles after the first division of the hyperrectangle. Since triangles are $2$-simplices, the polytope division thus only consists of simplices in this $2$-dimensional case. However, this is no longer the case in higher-dimensional settings.

\begin{figure}[h]
\centering
    \begin{subfigure}[b]{0.3\linewidth}
        \centering
        \input{Figures/S1_PDM-3d}
        \caption{}
        \label{fig: sf1: PDM3d}
    \end{subfigure}
    \hfill
    \begin{subfigure}[b]{0.3\linewidth}
        \centering
        \input{Figures/S2_PDM-3d}
        \caption{}
        \label{fig: sf2: PDM3d}
    \end{subfigure}
    \hfill
    \begin{subfigure}[b]{0.38\linewidth}
        \centering
        \input{Figures/S3_PDM-3d}
        \caption{}
        \label{fig: sf3: PDM3d}
    \end{subfigure}
    \hfill
\caption{Polytope division in $3\rm{d}$. (a) First polytope division with one polytope shaded. (b) Polytope division of one sub-polytope. (c) Dissection of the polytopes after the update step.}
\label{fig: PDM3d}
\end{figure}
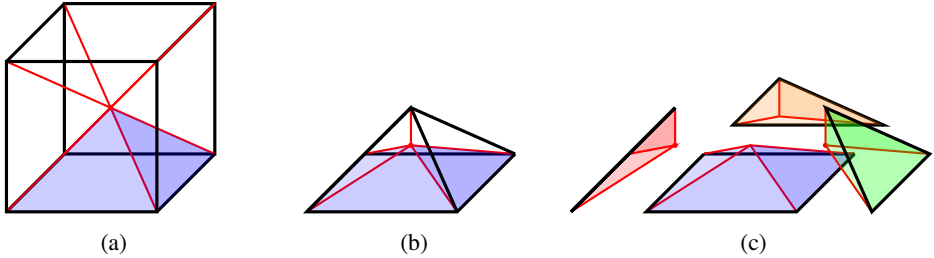

Figure~\ref{fig: PDM3d} displays the polytope division in the $3$-dimensional case based on facet linking after the first parameter selection and after one refinement step. Here, we distinguish two geometrical shapes: a pyramid with a square base and a simplex. Figure \ref{fig: sf2: PDM3d} displays the polytope division after one step. Figure \ref{fig: sf3: PDM3d} dissects the geometrical shapes in the polytope division. We omitted the frontal simplex here for the visibility of the other shapes. Note that the number of different types of geometrical shapes increases with increasing dimension.

As stated in the second point, determining the facets of a general polytope described by its vertices is non-trivial. However, this can still be done in PDM because we start with a hyperrectangle and know how to update the polytope division in every step. This leads to specific polytopes, which allows us to compute the facets quickly.

\begin{theorem}
\label{thm: main theorem polytope division}
    Let $P$ be a $d$-dimensional hyperrectangle and let $\calD^N$ be the polytope division of $P$ constructed by Algorithm~\ref{pdm} at step $N \ge 1$. Then
    $$\calD^N = \calS \cup \calB,$$
    where $\calS \subset \simp$ is a family of simplices and $\calB \subset \bp$ a family of $P$-boundary polytopes.
\end{theorem}
This observation is useful because we claim that we can determine the facets of both objects. The proof of this theorem is given at the end of this section and relies on the following intermediate results.

\medskip

For simplices, we have
\begin{lemma}
    \label{lemma: facets of simplices}
    Let $\sigma$ be a simplex with vertices $v_1, \ldots, v_{d+1}$. Then its facets are given by
    \begin{align}
    \partial \sigma = \bigl\{\conv\bigl(\{v_1, \ldots, v_{d+1}\} \backslash \{v_j\}\bigr)\, \big| \,j = 1, \ldots, d+1\bigr\}.
\end{align}
\end{lemma}

A similar result holds for $P$-boundary polytopes. We postpone its proof to Appendix~\ref{app: proofs}.
\begin{lemma}
\label{lemma: facets p-boundary polytopes}
    Let $d \ge 2$ and let $\tau \coloneq \conv\{\{x_1, \ldots, x_n\}\cup F\}$ be a $P$-boundary polytope with $F\in \partial_n P$ for some $1 \le n \le d-1$. Then,
    \begin{align*}
        \partial \tau = \left(\bigcup_{i = 1}^n \conv\bigl(\{x_1, \ldots, x_n\} \backslash \{x_i\} \cup F\bigr)\right)\cup \left(\bigcup_{G \in \partial F} \conv\bigl(\{x_1, \ldots, x_n\}\cup G\bigr)\right).
    \end{align*}
\end{lemma}

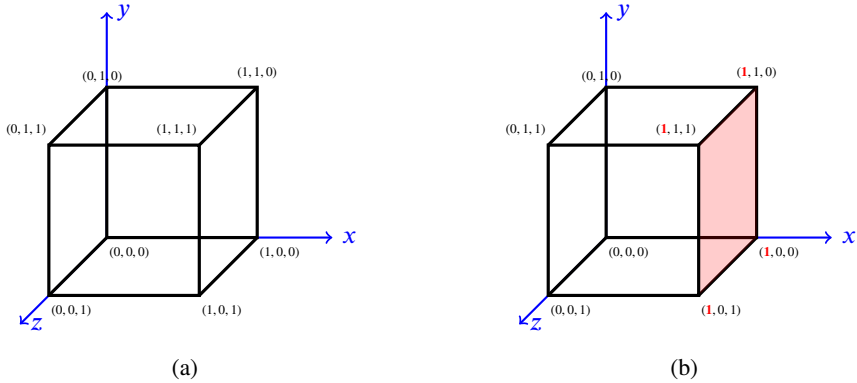
\begin{figure}[h]
\centering
    \begin{subfigure}[b]{0.49\linewidth}
        \centering
        \input{Figures/S1_Facets-example}
        \caption{}
        \label{fig: sf1: facets-example}
    \end{subfigure}
    \hfill
    \begin{subfigure}[b]{0.49\linewidth}
        \centering
        \input{Figures/S2_Facets-example}
        \caption{}
        \label{fig: sf2: facets-example}
    \end{subfigure}
\caption{Illustration of facets in hypercubes. (a) Example of a $3\rm{d}$-hypercube with vertices. (b)  Example of one facet of a 3-dimensional hypercube.}
\label{fig: facets-example}
\end{figure}

Lastly, to start Algorithm \ref{pdm}, we need to determine the facets of a hyperrectangle. This is relatively straightforward, as is implied in the following lemma.
\begin{lemma}
    \label{lemma: facets hyperrectangle}
    Let $P = [l_1, u_1]\times\cdots\times[l_d,u_d]$, $l_i$, $u_i \in \R$ be a $d$-dimensional hyper\-rectangle. Then,
    \begin{align*}
        \partial P =\bigcup_{i=1}^d \bigr\{F_{l_i}, F_{u_i}\bigl\}, 
    \end{align*}
    where for each $i \in \{1, \ldots, d\}$,
    \begin{align*}
        F_{l_i} &= \bigl\{p \in P 
        \; | \; \text{the $i$-th component of $p = l_i$}\bigr\},\\
        F_{u_i} &= \bigl\{p \in P\; |\; \text{the $i$-th component of $p = u_i$}\bigr\}.
    \end{align*}  
\end{lemma}
Figure \ref{fig: facets-example} provides a $3$-dimensional example for the hypercube $[0,1]^3$ with facet $F_{u_1}$ marked.

The motivation for the facet linking procedure is the scalability of the number of facets in higher dimensions. We know that every $d$-dimensional hyperrectangle has $2d$ facets. Since we only consider simplices and $P$-boundary polytopes, we can show that this beneficial scaling property still holds true in later steps of Algorithm~\ref{pdm}.

\begin{corollary}
    Let $P$ be a $d$-dimensional hyperrectangle, let $\calD^N$ be the polytope division of $P$ at step $N\ge 1$ of Algorithm~\ref{pdm}, and let $D \in \calD^N$. It then follows that $|\partial D| \le 2d$, where $|\cdot|$ denotes the cardinality of a set.
\end{corollary}

\begin{proof}
    $D$ is either a simplex or a $P$-boundary polytope.
    If $D$ is a simplex, it has $d+1 \le 2d$ facets.
    If $D$ is a $P$-boundary polytope, then it has the form $\conv(\{x_1, \ldots, x_n\} \cup F)$ for some $n \in \N$, $1 \le n \le d-1$, where $F \in \partial_n P$. Since $\left|\partial F\right| = 2(d-n)$, we have by Lemma~\ref{lemma: facets p-boundary polytopes} that
    \begin{align*}
        \left|\partial D \right| \le n + 2(d-n) \le 2d. 
    \end{align*}
\end{proof}

\begin{figure}[h]
\centering
    \begin{subfigure}[b]{0.32\linewidth}
        \centering
        \input{Figures/ConvexHull1}
        \caption{}
        \label{fig: sf1: convexhull}
    \end{subfigure}
    \hfill
    \begin{subfigure}[b]{0.32\linewidth}
        \centering
        \input{Figures/ConvexHull2}
        \caption{}
        \label{fig: sf2: convexhull}
    \end{subfigure}
    \hfill
    \begin{subfigure}[b]{0.32\linewidth}
        \centering
        \input{Figures/ConvexHull3}
        \caption{}
        \label{fig: sf3: convexhull}
    \end{subfigure}
\caption{A 2-\rm{d} convex hull example with point $v_5$ removed (a) $\conv(\{\sfb, v_1, \ldots, v_5\} \backslash \{v_5\})$; (b) $\conv(\{v_1, \ldots, v_5\} \backslash \{v_5\})$; (c) $\conv\big(\sfb \cup \conv(\{v_1, \ldots, v_5\} \backslash \{v_5\})\big)$. Note that this example never occurs as an outcome of Algorithm~\ref{pdm}, since no polytope in a $2$-dimensional polytope division in PDM ever has $5$ vertices.}
\label{fig: convexhull}
\end{figure}
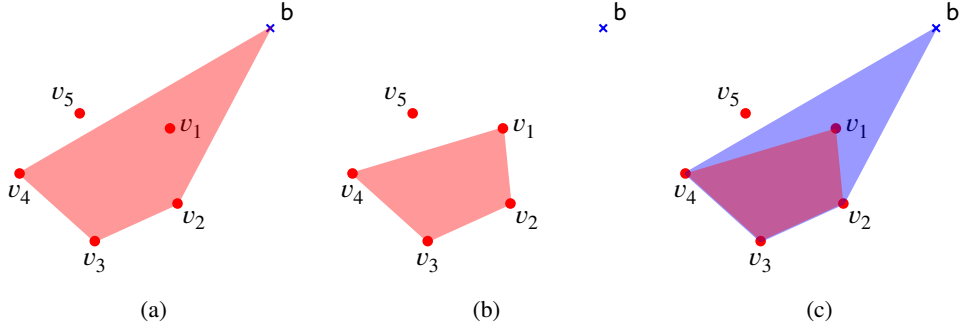

With these intermediate results, we can now prove Theorem \ref{thm: main theorem polytope division}.
\begin{proof}[Proof of Theorem~\ref{thm: main theorem polytope division}]
    We prove this by induction over $N$. Let $N=1$, hence we look at the first step of PDM. Via the facet linking operations, it holds that
    \begin{align*}
        \calD^1 = \bigcup_{F \in \partial P} \conv(p \cup F),\qquad p\in \interior{P}.
    \end{align*}
    If $d \le 2$, then $\calD^1 \subset \simp$ (i.e., $\calD^1$ is a family of simplices), while for $d \ge 2$, $\calD^1 \subset \bp$ (i.e., $\calD^1$ is a family of $P$-boundary polytopes).
    Suppose for $N \ge 1$, it holds that $\calD^N = \calS \cup \calB$, with $\calS \subset \simp$ and $\calB \subset \bp$, then we show that it also holds for $N + 1$.
    In one step of Algorithm~\ref{pdm}, one of the polytopes in $\calD^N$ is replaced via the facet linking procedure.
    We now consider the two cases. In the first case, suppose we replace a simplex $\sigma \in \calS$ having vertices $v_1,\ldots,v_{d+1}$. Then we know that \begin{align*}
    \FL(\sfb_\sigma, \sigma) &= \left\{\conv\big(\sfb_\sigma \cup \conv(\{v_1, \ldots, v_{d+1}\}\backslash \{v_j\})\big) \; \big | \; j = 1, \ldots, d+ 1\right\} \\
    &= \left\{\conv(\{\sfb_\sigma, v_1, \ldots, v_{d+1}\} \backslash \{v_j
    \}) \; \big| \; j = 1, \ldots, d+1\right\} \subset \simp.
    \end{align*}
    This statement is true by the definition of the convex hull (see Figure~\ref{fig: convexhull} for an illustration). Hence, if we refine a simplex, we know that the updated set $\calD^{N+1}$ is again a union of simplices and $P$-boundary polytopes.
    In the second case, suppose we replace a $P$-boundary polytope $\tau$ in step $N$. Then, we know $\tau \coloneq \conv(\{x_1, \ldots, x_n\} \cup F\})$ for some $F\in \partial_n P$ and, from Lemma~\ref{lemma: facets p-boundary polytopes}, the facets have either the form 
    \[
        \text{$\conv(\{x_1, \ldots, x_n\} \backslash \{x_i\} \cup F)$ or $\conv(\{x_1, \ldots, x_n\} \cup G)$ with $G \in \partial F$.}
    \]
    Consequently, the new polytopes take either the form
    \begin{align*}
        \conv\big(\sfb_\tau \cup \conv(\{x_1, \ldots, x_n\} \backslash \{x_i\}, F)\big) = \conv(\{\sfb_\tau, x_1, \ldots, x_n\} \backslash \{x_i\} \cup F) \subset \bp,
    \end{align*}
    or
    \begin{align*}
        \conv(\sfb_\tau \cup \conv(\{x_1, \ldots, x_n\} \cup G)) = \conv(\{\sfb_\tau, x_1, \ldots, x_n\} \cup G),  
    \end{align*}
    where $G \in \partial F \subset \partial_{n+1} P$, for $1 \le n \le d-1$. If $n = d-1$ then $G$ is a vertex and therefore $\conv(\{\sfb_\tau, x_1, \ldots, x_n\} \cup G) \in \simp$. If $n < d-1$, then $\conv(\{\sfb_\tau, x_1, \ldots, x_n\} \cup G) \in \bp$.
    This concludes the proof.
\end{proof}

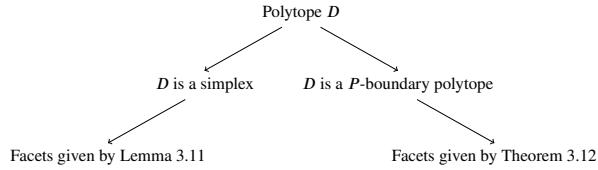
\begin{figure}[h]
\centering
\resizebox{8cm}{!}{\input{Figures/FacetSchedule}}
\caption{Schedule to determine the facets}\label{fig: facet schedule}
\end{figure}

Figure \ref{fig: facet schedule} presents a short overview of the procedure to determine the facets.

%% file: Figures/S1_PDM-2d.tex
\begin{tikzpicture}
        \def\d{2};
        \def\o{2};

        \coordinate (S1bl) at (0,0);
        \coordinate (S1br) at (\d,0);
        \coordinate (S1tr) at (\d,\d);
        \coordinate (S1tl) at (0,\d);
        \coordinate (S1c) at (\d/2,\d/2);
        \coordinate (bc1) at (\d/4, \d/2);
        \coordinate (bc2) at (\d/2, \d/4);
        \coordinate (bc3) at (3*\d/4, \d/2);
        \coordinate (bc4) at (\d/2, 3*\d/4);

        \draw[thick, color=black] (S1bl) -- (S1br) -- (S1tr) -- (S1tl) -- (S1bl);
        \draw[thick, color=red] (S1bl) -- (S1c) -- (S1br);
        \draw[thick, color=red] (S1tl) -- (S1c) -- (S1tr);
        \draw[draw=none, fill=red] (S1c) circle (2pt);

        \coordinate (S1trc) at (3*\d/4,3*\d/4);
        \coordinate (S1tc) at (\d/2,\d);
        \coordinate (S1tlc) at (\d/4,3*\d/4);
\end{tikzpicture}

%% file: Figures/S2_PDM-2d.tex
\begin{tikzpicture}
        \def\d{2};
        \def\o{2};

        \coordinate (S1bl) at (0,0);
        \coordinate (S1br) at (\d,0);
        \coordinate (S1tr) at (\d,\d);
        \coordinate (S1tl) at (0,\d);
        \coordinate (S1c) at (\d/2,\d/2);
        \coordinate (bc1) at (\d/6, \d/2);
        \coordinate (bc2) at (\d/2, \d/6);
        \coordinate (bc3) at (5*\d/6, \d/2);
        \coordinate (bc4) at (\d/2, 5*\d/6);

        \draw[draw=none, fill = blue] plot[mark=x] coordinates{(bc1)};
        \draw[draw=none, fill = blue] plot[mark=x] coordinates{(bc2)};
        \draw[draw=none, fill = blue] plot[mark=x] coordinates{(bc3)};
        \draw[draw=none, fill = blue] plot[mark=x] coordinates{(bc4)};

        \draw[thick, color=black] (S1bl) -- (S1br) -- (S1tr) -- (S1tl) -- (S1bl);
        \draw[thick, color=red] (S1bl) -- (S1c) -- (S1br);
        \draw[thick, color=red] (S1tl) -- (S1c) -- (S1tr);
        \draw[draw=none, fill=red] (S1c) circle (2pt);

        \coordinate (S1trc) at (3*\d/4,3*\d/4);
        \coordinate (S1tc) at (\d/2,\d);
        \coordinate (S1tlc) at (\d/4,3*\d/4);

\end{tikzpicture}

%% file: Figures/S3_PDM-2d.tex
\begin{tikzpicture}
        \def\d{2};
        \def\o{2};

        \coordinate (S1bl) at (0,0);
        \coordinate (S1br) at (\d,0);
        \coordinate (S1tr) at (\d,\d);
        \coordinate (S1tl) at (0,\d);
        \coordinate (S1c) at (\d/2,\d/2);
        \coordinate (bc1) at (\d/6, \d/2);
        \coordinate (bc2) at (\d/2, \d/6);
        \coordinate (bc3) at (5*\d/6, \d/2);
        \coordinate (bc4) at (\d/2, 5*\d/6);

        \draw[draw=none, fill = blue] plot[mark=x] coordinates{(bc1)};
        \draw[draw=none, fill = blue] plot[mark=x] coordinates{(bc2)};
        \draw[draw=none, fill = blue] plot[mark=x] coordinates{(bc3)};
        \draw[draw=none, fill = blue] plot[mark=x] coordinates{(bc4)};

        \draw[thick, color=black] (S1bl) -- (S1br) -- (S1tr) -- (S1tl) -- (S1bl);
        \draw[thick, color=red] (S1bl) -- (S1c) -- (S1br);
        \draw[thick, color=red] (S1tl) -- (S1c) -- (S1tr);
        \draw[draw=none, fill=red] (S1c) circle (2pt);
        \draw[draw=blue, fill= none] (bc1) circle (4pt);

        \coordinate (S1trc) at (3*\d/4,3*\d/4);
        \coordinate (S1tc) at (\d/2,\d);
        \coordinate (S1tlc) at (\d/4,3*\d/4);
\end{tikzpicture}

%% file: Figures/S4_PDM-2d.tex
\begin{tikzpicture}
        \def\d{2};
        \def\o{2};

        \coordinate (S1bl) at (0,0);
        \coordinate (S1br) at (\d,0);
        \coordinate (S1tr) at (\d,\d);
        \coordinate (S1tl) at (0,\d);
        \coordinate (S1c) at (\d/2,\d/2);
        \coordinate (bc1) at (\d/6, \d/2);
        \coordinate (bc2) at (\d/2, \d/6);
        \coordinate (bc3) at (5*\d/6, \d/2);
        \coordinate (bc4) at (\d/2, 5*\d/6);

        \draw[thick, color=black] (S1bl) -- (S1br) -- (S1tr) -- (S1tl) -- (S1bl);
        \draw[thick, color=red] (S1bl) -- (S1c) -- (S1br);
        \draw[thick, color=red] (S1tl) -- (S1c) -- (S1tr);
        \draw[thick, color=red] (S1bl) -- (bc1) -- (S1c);
        \draw[thick, color=red] (bc1) -- (S1tl);
        \draw[draw=none, fill=red] (S1c) circle (1pt);
        \draw[draw=none, fill=red] (bc1) circle (1pt);

        \coordinate (S1trc) at (3*\d/4,3*\d/4);
        \coordinate (S1tc) at (\d/2,\d);
        \coordinate (S1tlc) at (\d/4,3*\d/4);
\end{tikzpicture}

%% file: Figures/S5_PDM-2d.tex
\begin{tikzpicture}
        \def\d{2};
        \def\o{2};

        \coordinate (S1bl) at (0,0);
        \coordinate (S1br) at (\d,0);
        \coordinate (S1tr) at (\d,\d);
        \coordinate (S1tl) at (0,\d);
        \coordinate (S1c) at (\d/2,\d/2);
        \coordinate (bc1) at (\d/6, \d/2);
        \coordinate (bc2) at (\d/2, \d/6);
        \coordinate (bc3) at (5*\d/6, \d/2);
        \coordinate (bc4) at (\d/2, 5*\d/6);

        \draw[thick, color=black] (S1bl) -- (S1br) -- (S1tr) -- (S1tl) -- (S1bl);
        \draw[thick, color=red] (S1bl) -- (S1c) -- (S1br);
        \draw[thick, color=red] (S1tl) -- (S1c) -- (S1tr);
        \draw[thick, color=red] (S1bl) -- (bc1) -- (S1c);
        \draw[thick, color=red] (bc1) -- (S1tl);
        \draw[draw=none, fill=red] (S1c) circle (1pt);
        \draw[draw=none, fill=red] (bc1) circle (1pt);

        \coordinate (bc11) at (\d/18,\d/2);
        \coordinate (bc12) at (\d/4,\d/3);
        \coordinate (bc13) at (\d/4,2*\d/3);

        \draw[draw=black, fill=none] plot[mark=x] coordinates{(bc11)};
        \draw[draw=black, fill=none] plot[mark=x] coordinates{(bc12)};
        \draw[draw=black, fill=none] plot[mark=x] coordinates{(bc13)};
        \draw[draw=black, fill=none] plot[mark=x] coordinates{(bc2)};
        \draw[draw=black, fill=none] plot[mark=x] coordinates{(bc3)};
        \draw[draw=black, fill=none] plot[mark=x] coordinates{(bc4)};
\end{tikzpicture}

%% file: Figures/P.tex
\begin{tikzpicture}
        \def\d{2};
        \def\o{2};

        \coordinate (S1bl) at (0,0);
        \coordinate (S1br) at (\d,0);
        \coordinate (S1tr) at (\d,\d);
        \coordinate (S1tl) at (0,\d);

        \draw[draw=none, fill=black] (S1bl) circle (1pt);
        \draw[draw=none, fill=black] (S1br) circle (1pt);
        \draw[draw=none, fill=black] (S1tr) circle (1pt);
        \draw[draw=none, fill=black] (S1tl) circle (1pt);

        \draw[thick, color=black] (S1bl) -- (S1br) -- (S1tr) -- (S1tl) -- (S1bl);

\end{tikzpicture}

%% file: Figures/P_facets.tex
\begin{tikzpicture}
        \def\d{2};
        \def\o{2};

        \coordinate (S1bl) at (0,0);
        \coordinate (S1br) at (\d,0);
        \coordinate (S1tr) at (\d,\d);
        \coordinate (S1tl) at (0,\d);

        \coordinate (f1) at (0, \d/2);
        \coordinate (f2) at (\d, \d/2);
        \coordinate (f3) at (\d/2, \d);
        \coordinate (f4) at (\d/2, 0);

        \draw[draw=none, fill=black] (S1bl) circle (1pt);
        \draw[draw=none, fill=black] (S1br) circle (1pt);
        \draw[draw=none, fill=black] (S1tr) circle (1pt);
        \draw[draw=none, fill=black] (S1tl) circle (1pt);

        \draw[thick, color=black] (S1bl) -- (S1br) -- (S1tr) -- (S1tl) -- (S1bl);
        \draw[draw=red, fill=none] (f1) ellipse (3pt and 32pt);
        \draw[draw=red, fill=none] (f2) ellipse (3pt and 32pt);
        \draw[draw=red, fill=none] (f3) ellipse (32pt and 3pt);
        \draw[draw=red, fill=none] (f4) ellipse (32pt and 3pt);                       
\end{tikzpicture}

%% file: Figures/P_vertices.tex
\begin{tikzpicture}
        \def\d{2};
        \def\o{2};

        \coordinate (S1bl) at (0,0);
        \coordinate (S1br) at (\d,0);
        \coordinate (S1tr) at (\d,\d);
        \coordinate (S1tl) at (0,\d);

        \draw[draw=none, fill=black] (S1bl) circle (1pt);
        \draw[draw=none, fill=black] (S1br) circle (1pt);
        \draw[draw=none, fill=black] (S1tr) circle (1pt);
        \draw[draw=none, fill=black] (S1tl) circle (1pt);

        \draw[thick, color=black] (S1bl) -- (S1br) -- (S1tr) -- (S1tl) -- (S1bl);
        \draw[draw=red, fill=none] (S1bl) circle (4pt);
        \draw[draw=red, fill=none] (S1br) circle (4pt);
        \draw[draw=red, fill=none] (S1tl) circle (4pt);
        \draw[draw=red, fill=none] (S1tr) circle (4pt);                              
\end{tikzpicture}

%% file: Figures/FS1point-facet.tex
\begin{tikzpicture}
        \def\d{4};
        \def\o{2};

        \coordinate (S1bl) at (0,0);
        \coordinate (S1br) at (\d,0);
        \coordinate (S1tr) at (\d,\d);
        \coordinate (S1tl) at (0,\d);
        \coordinate (S1c) at (\d/2,\d/2);
        \coordinate (p) at (\d/2 + \d/18, \d/2);
        \coordinate (F) at (\d/2, \d/18);
        \coordinate (bc1) at (\d/4, \d/2);
        \coordinate (bc2) at (\d/2, \d/4);
        \coordinate (bc3) at (3*\d/4, \d/2);
        \coordinate (bc4) at (\d/2, 3*\d/4);

        \draw[dashed, color=black] (S1bl) -- (S1br) -- (S1tr) -- (S1tl) -- (S1bl);
        \draw[very thick, color=red] (S1bl) -- (S1br);
        \draw[draw=none, fill=red] (S1c) circle (2pt);
        \node at (p) {$p$};
        \node at (F) {$F$};
\end{tikzpicture}

%% file: Figures/FS2convex-hull.tex
\begin{tikzpicture}
        \def\d{4};
        \def\o{2};

        \coordinate (S1bl) at (0,0);
        \coordinate (S1br) at (\d,0);
        \coordinate (S1tr) at (\d,\d);
        \coordinate (S1tl) at (0,\d);
        \coordinate (S1c) at (\d/2,\d/2);
        \coordinate (bc1) at (\d/4, \d/2);
        \coordinate (bc2) at (\d/2, \d/4);
        \coordinate (bc3) at (3*\d/4, \d/2);
        \coordinate (bc4) at (\d/2, 3*\d/4);

        \draw[dashed, color=black] (S1bl) -- (S1br) -- (S1tr) -- (S1tl) -- (S1bl);

        \draw[very thick, color=red] (S1bl) -- (S1br);
        \draw[draw=none, fill=red] (S1c) circle (2pt);
        \draw[draw = none, fill = red, opacity=0.4] (S1bl) -- (S1c) -- (S1br);
        \draw[very thin, color=red] (S1bl) -- (S1c) -- (S1br);
\end{tikzpicture}

%% file: Figures/S1_PDM-3d.tex
\begin{tikzpicture}
		[cube/.style={very thick,black},
			grid/.style={very thin,gray},
			axis/.style={->,blue,thick}]

   \coordinate (lbb) at (0,0,0);
   \coordinate (lbf) at (0,0,2);
   \coordinate (rbb) at (2,0,0);
   \coordinate (rbf) at (2,0,2);
   \coordinate (ltb) at (0,2,0);
   \coordinate (ltf) at (0,2,2);
   \coordinate (rtb) at (2,2,0);
   \coordinate (rtf) at (2,2,2);
   \coordinate (c1) at (1,1,1);

	\draw[cube] (lbb) -- (ltb) -- (rtb) -- (rbb) -- cycle;

    \draw[cube] (lbb) -- (lbf);
    \draw[cube] (rtb) -- (rtf);
    \draw[thick, color=red] (ltf) -- (c1) -- (rtf);
    \draw[thick, color=red] (ltb) -- (c1) -- (rtb);
    \draw[thick, color=red] (lbf) -- (c1) -- (lbb);
    \draw[thick, color=red] (lbb) -- (c1) -- (rbb);
    \draw[thick, color=red] (rbf) -- (c1);

	\draw[cube] (ltb) -- (ltf);
	\draw[cube] (rbb) -- (rbf);
    \draw[draw=none, fill=blue, opacity=.2] (lbf) -- (c1) -- (rbf);
    \draw[draw=none, fill=blue, opacity=.3] (rbf) -- (c1) -- (rbb);
    \draw[draw=none, fill=red] (c1) circle (1pt);
    \draw[cube] (lbf) -- (ltf) -- (rtf) -- (rbf) -- cycle;
\end{tikzpicture}

%% file: Figures/S2_PDM-3d.tex
\begin{tikzpicture}
		[cube/.style={very thick,black},
			grid/.style={very thin,gray},
			axis/.style={->,blue,thick}]

   \coordinate (lbb) at (0,0,0);
   \coordinate (lbf) at (0,0,2);
   \coordinate (rbb) at (2,0,0);
   \coordinate (rbf) at (2,0,2);
   \coordinate (ltb) at (0,2,0);
   \coordinate (ltf) at (0,2,2);
   \coordinate (rtb) at (2,2,0);
   \coordinate (rtf) at (2,2,2);
   \coordinate (c1) at (1,1,1);
   \coordinate (c2) at (1,1/2,1);

    \draw[cube] (rbb) -- (lbb);
    \draw[thick, color=red] (c1) -- (c2);
    \draw[thick, color=red] (lbb) -- (c2) -- (rbb);
    \draw[thick, color=red] (lbf) -- (c2) -- (rbf);

    \draw[draw=none, fill=blue, opacity=.15] (lbb) -- (c2) -- (lbf);
    \draw[draw=none, fill=blue, opacity=.2] (lbf) -- (c2) -- (rbf);
    \draw[draw=none, fill=blue, opacity=.3] (rbf) -- (c2) -- (rbb);

    \draw[cube] (lbb) -- (lbf) -- (rbf) -- (rbb);
    \draw[cube] (lbf) -- (c1) -- (rbf);
    \draw[cube] (lbb) -- (c1) -- (rbb);
    \draw[draw=none, fill=red] (c2) circle (1pt);
\end{tikzpicture}

%% file: Figures/S3_PDM-3d.tex
\begin{tikzpicture}
		[cube/.style={very thick,black},
			grid/.style={very thin,gray},
			axis/.style={->,blue,thick}]

   \coordinate (lbb) at (0,0,0);
   \coordinate (lbf) at (0,0,2);
   \coordinate (rbb) at (2,0,0);
   \coordinate (rbf) at (2,0,2);
   \coordinate (ltb) at (0,2,0);
   \coordinate (ltf) at (0,2,2);
   \coordinate (rtb) at (2,2,0);
   \coordinate (rtf) at (2,2,2);
   \coordinate (c1) at (1,1,1);
   \coordinate (c2) at (1,1/2,1);

   \coordinate (s1lbb) at (-1,0,0);
   \coordinate (s1lbf) at (-1,0,2);
   \coordinate (s1c1) at (0,1,1);
   \coordinate (s1c2) at (0,1/2,1);

   \coordinate (s2rbb) at (3,0,0);
   \coordinate (s2rbf) at (3,0,2);
   \coordinate (s2c1) at (2,1,1);
   \coordinate (s2c2) at (2,1/2,1);

   \coordinate (s3lbf) at (-1/2,3/2,2);
   \coordinate (s3rbf) at (3/2,3/2,2);
   \coordinate (s3c1) at (1/2,5/2,1);
   \coordinate (s3c2) at (1/2,2,1);

   \coordinate (s4lbb) at (0,0,-1);
   \coordinate (s4rbb) at (2,0,-1);
   \coordinate (s4c1) at (1,1,0);
   \coordinate (s4c2) at (1, 1/2, 0);

    \draw[cube] (rbb) -- (lbb);
    \draw[thick, color=red] (lbb) -- (c2) -- (rbb);
    \draw[thick, color=red] (lbf) -- (c2) -- (rbf);

    \draw[draw=none, fill=blue, opacity=.15] (lbb) -- (c2) -- (lbf);
    \draw[draw=none, fill=blue, opacity=.2] (lbf) -- (c2) -- (rbf);
    \draw[draw=none, fill=blue, opacity=.3] (rbf) -- (c2) -- (rbb);

    \draw[cube] (lbb) -- (lbf) -- (rbf) -- (rbb);
  
    \draw[thick, color=red] (s1lbb) -- (s1c2) -- (s1lbf);
    \draw[thick, color=red] (s1c1) -- (s1c2);
    \draw[draw=none, fill=red] (s1c2) circle (1pt);
    \draw[draw=none, fill=red, opacity=.15] (s1lbb) -- (s1c2) -- (s1c1);
    \draw[draw=none, fill=red, opacity=.2] (s1lbf) -- (s1c2) -- (s1c1);
    \draw[cube] (s1lbb) -- (s1lbf) -- (s1c1) -- cycle;

    \draw[thick, color=red] (s4lbb) -- (s4c2) -- (s4rbb);
    \draw[thick, color=red] (s4c1) -- (s4c2);
    \draw[cube] (s4lbb) -- (s4rbb) -- (s4c1) -- cycle;
    \draw[draw=none, fill=orange, opacity=.15] (s4lbb) -- (s4c2) -- (s4rbb);
    \draw[draw=none, fill=orange, opacity=.2] (s4lbb) -- (s4c2) -- (s4c1);
    \draw[draw=none, fill=orange, opacity=.3] (s4rbb) -- (s4c1) -- (s4c2);

    \draw[thick, color=red] (s2rbb) -- (s2c2) -- (s2rbf);
    \draw[thick, color=red] (s2c1) -- (s2c2);
    \draw[draw=none, fill=red] (s2c2) circle (1pt);
    \draw[draw=none, fill=green, opacity=.15] (s2rbb) -- (s2c1) -- (s2c2);
    \draw[draw=none, fill=green, opacity=.2] (s2rbf) -- (s2c1) -- (s2c2);
    \draw[draw=none, fill=green, opacity=.3] (s2rbb) -- (s2c1) -- (s2rbf);
    \draw[cube] (s2rbb) -- (s2rbf) -- (s2c1) -- cycle;
\end{tikzpicture}

%% file: Figures/S1_Facets-example.tex
\begin{tikzpicture}
		[cube/.style={very thick,black},
			grid/.style={very thin,gray},
			axis/.style={->,blue,thick}]

   \coordinate (lbb) at (0,0,0);
   \coordinate (lbf) at (0,0,2);
   \coordinate (rbb) at (2,0,0);
   \coordinate (rbf) at (2,0,2);
   \coordinate (ltb) at (0,2,0);
   \coordinate (ltf) at (0,2,2);
   \coordinate (rtb) at (2,2,0);
   \coordinate (rtf) at (2,2,2);

   \coordinate (Tlbb) at (0.3,-0.2,0);
   \coordinate (Tlbf) at (0.3,-0.2,2);
   \coordinate (Trbb) at (2.3,-0.2,0);
   \coordinate (Trbf) at (2.3,-0.2,2);
   \coordinate (Tltb) at (0,2.2,0.15);
   \coordinate (Tltf) at (-0.3,2.2,2);
   \coordinate (Trtb) at (2,2.2,0);
   \coordinate (Trtf) at (1.7,2.2,2);

	\draw[axis] (0,0,0) -- (3,0,0) node[anchor=west]{$x$};
	\draw[axis] (0,0,0) -- (0,3,0) node[anchor=west]{$y$};
	\draw[axis] (0,0,0) -- (0,0,3) node[anchor=west]{$z$};

	\draw[cube] (lbb) -- (ltb) -- (rtb) -- (rbb) -- cycle;

    \draw[cube] (lbb) -- (lbf);
    \draw[cube] (rtb) -- (rtf);

	\draw[cube] (ltb) -- (ltf);
	\draw[cube] (rbb) -- (rbf);
    
    \draw[cube] (lbf) -- (ltf) -- (rtf) -- (rbf) -- cycle;
   
    \node at (Tlbb) {\tiny $(0,0,0)$};
    \node at (Tlbf) {\tiny $(0,0,1)$};
    \node at (Trbb) {\tiny $(1,0,0)$};
    \node at (Trbf) {\tiny $(1,0,1)$};
    \node at (Tltb) {\tiny $(0,1,0)$};
    \node at (Tltf) {\tiny $(0,1,1)$};
    \node at (Trtb) {\tiny $(1,1,0)$};
    \node at (Trtf) {\tiny $(1,1,1)$};
    
\end{tikzpicture}

%% file: Figures/S2_Facets-example.tex
\begin{tikzpicture}
		[cube/.style={very thick,black},
			grid/.style={very thin,gray},
			axis/.style={->,blue,thick}]

   \coordinate (lbb) at (0,0,0);
   \coordinate (lbf) at (0,0,2);
   \coordinate (rbb) at (2,0,0);
   \coordinate (rbf) at (2,0,2);
   \coordinate (ltb) at (0,2,0);
   \coordinate (ltf) at (0,2,2);
   \coordinate (rtb) at (2,2,0);
   \coordinate (rtf) at (2,2,2);

   \coordinate (Tlbb) at (0.3,-0.2,0);
   \coordinate (Tlbf) at (0.3,-0.2,2);
   \coordinate (Trbb) at (2.3,-0.2,0);
   \coordinate (Trbf) at (2.3,-0.2,2);
   \coordinate (Tltb) at (0,2.2,0.15);
   \coordinate (Tltf) at (-0.3,2.2,2);
   \coordinate (Trtb) at (2,2.2,0);
   \coordinate (Trtf) at (1.7,2.2,2);
   
	\draw[axis] (0,0,0) -- (3,0,0) node[anchor=west]{$x$};
	\draw[axis] (0,0,0) -- (0,3,0) node[anchor=west]{$y$};
	\draw[axis] (0,0,0) -- (0,0,3) node[anchor=west]{$z$};

	\draw[cube] (lbb) -- (ltb) -- (rtb) -- (rbb) -- cycle;

    \draw[cube] (lbb) -- (lbf);
    \draw[cube] (rtb) -- (rtf);

	\draw[cube] (ltb) -- (ltf);
	\draw[cube] (rbb) -- (rbf);
 
    \draw[draw=none, fill=red, opacity=.2] (rbf) -- (rbb) -- (rtb) -- (rtf);
    
    \draw[cube] (lbf) -- (ltf) -- (rtf) -- (rbf) -- cycle;
   
    \node at (Tlbb) {\tiny $(0,0,0)$};
    \node at (Tlbf) {\tiny $(0,0,1)$};
    \node at (Trbb) {\tiny $(\textcolor{red}{\mathbf{1}},0,0)$};
    \node at (Trbf) {\tiny $(\textcolor{red}{\mathbf{1}},0,1)$};
    \node at (Tltb) {\tiny $(0,1,0)$};
    \node at (Tltf) {\tiny $(0,1,1)$};
    \node at (Trtb) {\tiny $(\textcolor{red}{\mathbf{1}},1,0)$};
    \node at (Trtf) {\tiny $(\textcolor{red}{\mathbf{1}},1,1)$};
    
\end{tikzpicture}

%% file: Figures/ConvexHull1.tex
\begin{tikzpicture}
        \def\d{4};
        \def\o{2};

        \coordinate (S1bl) at (0,0);
        \coordinate (S1br) at (\d,0);
        \coordinate (S1tr) at (\d,\d);
        \coordinate (S1tl) at (0,\d);
        \coordinate (v1) at (\d/2,\d/2);
        \coordinate (v1t) at (\d/2 + \d/14, \d/2);
        \coordinate (v2) at (\d/2 + 1/10, \d/4);
        \coordinate (v2t) at (\d/2 + 1/10 + \d/18, \d/4 - \d/18);
        \coordinate (v3) at (\d/4,\d/8);
        \coordinate (v3t) at (\d/4, \d/8 - \d/14);
        \coordinate (v4) at (0, 3*\d/8 -1/10);
        \coordinate (v4t) at (0, 3*\d/8 -1/10 - \d/14);
        \coordinate (v5) at (\d/4 - 1/5,\d/2 + 1/5);
        \coordinate (v5t) at (\d/4 - 1/5 - \d/18,\d/2 + 1/5 + \d/18);
        \coordinate (bar) at (5*\d/6, 5*\d/6);
        \coordinate (tbar) at (5*\d/6 + \d/18, 5*\d/6 + \d/18);
        \coordinate (bc1) at (\d/4, \d/2);
        \coordinate (bc2) at (\d/2, \d/4);
        \coordinate (bc3) at (3*\d/4, \d/2);
        \coordinate (bc4) at (\d/2, 3*\d/4);

        \draw[draw=none, fill=red] (v1) circle (2pt);
        \draw[draw=none, fill=red]
        (v2) circle (2pt);
        \draw[draw=none, fill=red] (v3) circle (2pt);
        \draw[draw=none, fill=red]
        (v4) circle (2pt);
        \draw[draw=none, fill=red] (v5) circle (2pt);
        \draw[thick, draw=blue, fill = blue] plot[mark=x] coordinates{(bar)};

        \node at (v1t) {$v_1$};
        \node at (v2t) {$v_2$};
        \node at (v3t) {$v_3$};
        \node at (v4t) {$v_4$};
        \node at (v5t) {$v_5$};
        \node at (tbar) {$\sfb$};

        \draw[thick, draw = red, fill = red, opacity=0.4] (v4) -- (bar) -- (v2) -- (v3) -- (v4);
\end{tikzpicture}

%% file: Figures/ConvexHull2.tex
\begin{tikzpicture}
        \def\d{4};
        \def\o{2};

        \coordinate (S1bl) at (0,0);
        \coordinate (S1br) at (\d,0);
        \coordinate (S1tr) at (\d,\d);
        \coordinate (S1tl) at (0,\d);
        \coordinate (v1) at (\d/2,\d/2);
        \coordinate (v1t) at (\d/2 + \d/14, \d/2);
        \coordinate (v2) at (\d/2 + 1/10, \d/4);
        \coordinate (v2t) at (\d/2 + 1/10 + \d/18, \d/4 - \d/18);
        \coordinate (v3) at (\d/4,\d/8);
        \coordinate (v3t) at (\d/4, \d/8 - \d/14);
        \coordinate (v4) at (0, 3*\d/8 -1/10);
        \coordinate (v4t) at (0, 3*\d/8 -1/10 - \d/14);
        \coordinate (v5) at (\d/4 - 1/5,\d/2 + 1/5);
        \coordinate (v5t) at (\d/4 - 1/5 - \d/18,\d/2 + 1/5 + \d/18);
        \coordinate (bar) at (5*\d/6, 5*\d/6);
        \coordinate (tbar) at (5*\d/6 + \d/18, 5*\d/6 + \d/18);
        \coordinate (bc1) at (\d/4, \d/2);
        \coordinate (bc2) at (\d/2, \d/4);
        \coordinate (bc3) at (3*\d/4, \d/2);
        \coordinate (bc4) at (\d/2, 3*\d/4);

        \draw[draw=none, fill=red] (v1) circle (2pt);
        \draw[draw=none, fill=red]
        (v2) circle (2pt);
        \draw[draw=none, fill=red] (v3) circle (2pt);
        \draw[draw=none, fill=red]
        (v4) circle (2pt);
        \draw[draw=none, fill=red] (v5) circle (2pt);
        \draw[thick, draw=blue, fill = blue] plot[mark=x] coordinates{(bar)};

        \node at (v1t) {$v_1$};
        \node at (v2t) {$v_2$};
        \node at (v3t) {$v_3$};
        \node at (v4t) {$v_4$};
        \node at (v5t) {$v_5$};
        \node at (tbar) {$\sfb$};

        \draw[thick, draw = red, fill = red, opacity=0.4]  (v2) -- (v3) -- (v4) -- (v1) -- (v2);
\end{tikzpicture}

%% file: Figures/ConvexHull3.tex
\begin{tikzpicture}
        \def\d{4};
        \def\o{2};

        \coordinate (S1bl) at (0,0);
        \coordinate (S1br) at (\d,0);
        \coordinate (S1tr) at (\d,\d);
        \coordinate (S1tl) at (0,\d);
        \coordinate (v1) at (\d/2,\d/2);
        \coordinate (v1t) at (\d/2 + \d/14, \d/2);
        \coordinate (v2) at (\d/2 + 1/10, \d/4);
        \coordinate (v2t) at (\d/2 + 1/10 + \d/18, \d/4 - \d/18);
        \coordinate (v3) at (\d/4,\d/8);
        \coordinate (v3t) at (\d/4, \d/8 - \d/14);
        \coordinate (v4) at (0, 3*\d/8 -1/10);
        \coordinate (v4t) at (0, 3*\d/8 -1/10 - \d/14);
        \coordinate (v5) at (\d/4 - 1/5,\d/2 + 1/5);
        \coordinate (v5t) at (\d/4 - 1/5 - \d/18,\d/2 + 1/5 + \d/18);
        \coordinate (bar) at (5*\d/6, 5*\d/6);
        \coordinate (tbar) at (5*\d/6 + \d/18, 5*\d/6 + \d/18);
        \coordinate (bc1) at (\d/4, \d/2);
        \coordinate (bc2) at (\d/2, \d/4);
        \coordinate (bc3) at (3*\d/4, \d/2);
        \coordinate (bc4) at (\d/2, 3*\d/4);

        \draw[draw=none, fill=red] (v1) circle (2pt);
        \draw[draw=none, fill=red]
        (v2) circle (2pt);
        \draw[draw=none, fill=red] (v3) circle (2pt);
        \draw[draw=none, fill=red]
        (v4) circle (2pt);
        \draw[draw=none, fill=red] (v5) circle (2pt);
        \draw[thick, draw=blue, fill = blue] plot[mark=x] coordinates{(bar)};

        \node at (v1t) {$v_1$};
        \node at (v2t) {$v_2$};
        \node at (v3t) {$v_3$};
        \node at (v4t) {$v_4$};
        \node at (v5t) {$v_5$};
        \node at (tbar) {$\sfb$};

        \draw[thick, draw = blue, fill = blue, opacity=0.4] (v4) -- (bar) -- (v2) -- (v3) -- (v4);
        \draw[draw = none, fill = red, opacity=0.4] (v1) --  (v2) -- (v3) -- (v4) -- (v1);
\end{tikzpicture}

%% file: Figures/FacetSchedule.tex
\begin{tikzpicture}
        \def\d{2};

        \coordinate (T1) at (4,8);
        \coordinate (A11) at (3.6, 7.7);
        \coordinate (A12) at (2, 6.8);
        \coordinate (A21) at (4.4, 7.7);
        \coordinate (A22) at (6, 6.8);
        \coordinate (T2) at (2, 6.5);
        \coordinate (T3) at (6, 6.5);

        \coordinate (A31) at (1.6, 6.2);
        \coordinate (A32) at (0, 5.3);
        \coordinate (A41) at (6.4, 6.2);
        \coordinate (A42) at (8, 5.3);

        \coordinate (T4) at (0, 5);
        \coordinate (T5) at (8, 5);

        \coordinate (A51) at (6.1, 4.7);
        \coordinate (A52) at (4.4, 3.8);

        \coordinate (A61) at (6.9, 4.7);
        \coordinate (A62) at (8.6, 3.8);

        \coordinate (T6) at (4.4, 3.5);
        \coordinate (T7) at (8.6, 3.5);

        \coordinate (A71) at (4.4, 3.2);
        \coordinate (A72) at (6.1, 2.3);

        \coordinate (A81) at (8.6, 3.2);
        \coordinate (A82) at (6.9, 2.3);

        \coordinate (T8) at (6.5, 2);

        \node at (T1) {Polytope $D$};

        \draw[->] (A11) -- (A12);
        \draw[->] (A21) -- (A22);

        \node at (T2) {$D$ is a simplex};
        \node at (T3) {$D$ is a $P$-boundary polytope};

        \draw[->] (A31) -- (A32);
        \draw[->] (A41) -- (A42);

        \node at (T4) {Facets given by Lemma \ref{lemma: facets of simplices}};
        \node at (T5) {Facets given by Theorem \ref{lemma: facets p-boundary polytopes}};
\end{tikzpicture}

%% file: Reduced_Basis_Method.tex
The GSM and PDM can be used to construct a reduced basis. In this section, we provide a brief explanation of the reduced basis method. 

We consider a parametrized variational form
\begin{align}
\label{eq: parametrized variational form}
a(u(p), v; p) = f(v; p), \qquad \forall v \in \V,
\end{align}
where $\V$ is a Hilbert space, $a: \V \times \V \times P \to \R$ is a $\V$-coercive bilinear form and $f: \V \times P \to \R$ is a linear form. Here, $p$ is the parameter, and $\bar u(p)$ is the solution of \eqref{eq: parametrized variational form} depending on the parameter $p$. We assume the parameter set $P \subset \R^d$ is a $d$-dimensional hyperrectangle. In general, we know that $p \in P$, but the precise value of the parameter is unknown.

For a fixed parameter $p \in P$, the exact solution $\bar u(p)$ could be well approximated by the Finite Element Method (FEM). 
Based on the FEM, there exists a matrix $A(p) \in \R^{N_h \times N_h}$ and a vector $f(p) \in \R^{N_h}$ such that the FEM-solution $u(p)$ satisfies
\begin{align}\label{eq: FE-equation}
A(p) u(p) = f(p),
\end{align}
where $N_h\in\N$ is the number of degrees of freedom. Moreover, $A(p)$ is commonly referred to as the "stiffness" matrix, and $f(p)$ is the load vector (see e.g.\ \cite{zienkiewicz2005finite}).

In this setup, we assume that $A(p)$ and $f(p)$ have an affine decomposition. This means that there exist $Q_a, Q_f \in \N$ and scalar-valued functions $\{\theta_\eta^q(p)\,|\, q=1,\ldots,Q_\eta\}$, $\eta\in\{a,f\}$, such that $A(p)$ and $f(p)$ can be written as
\begin{align}
\label{ass: affine decomposition}
A(p) = \sum_{q=1}^{Q_a} \theta_a^q(p)A^q,\qquad
f(p) = \sum_{q=1}^{Q_f} \theta_f^{q}(p) f^{q}.
\end{align}
The benefit of this assumption is that each $A^q$ and $f^{q}$ can be pre-computed offline.

Although FEM could provide a fine approximate solution for known parameters $p \in P$, the typical high-dimensionality of the FE equation \eqref{eq: FE-equation} leads to very high computational costs, particularly in multi-query settings (see, e.g. \cite{benner2015survey} for a survey paper). To reduce the complexity of the solution space, we aim to construct a configuration $\bgamma \in \Gamma_N, \bgamma = \{p_1, \ldots, p_N\}$, leading to a \emph{reduced basis} $V_{\bgamma} = \textup{span}\{u^{p_1}, ..., u^{p_N}\}$. After orthonormalization of the elements of $V_{\bgamma}$, we also write $V_{\bgamma} = \textup{span}\{\xi_1, \ldots, \xi_N\}$. Then for every parameter $p \in P$, we assume there exist $\lambda_i, \,i \in \{1, \ldots, N\}$, such that $u(p) \approx \sum_{i=1}^N \lambda_i \, u^{p_i}$, where $N$ is the number of basis vectors in $V_{\bgamma}$ and, ideally, $N \ll N_h$.

A reduced basis $V_{\bgamma}$ can be beneficial because it allows us to construct the reduced system
\begin{align}
    \label{eq: reduced system}
    A_N(p) u_N(p) = f_N(p),
\end{align}
where $A_N(p) \coloneq V_{\bgamma}^\top A(p) V_{\bgamma} \in \R^{N \times N}$, $f_N(p) \coloneq V_{\bgamma}^\top f(p) \in \R^N$ and $u_N(p) \in \R^N$, thus allowing us to approximate $u(p) \approx u_{\bgamma}(p) \coloneq V_{\bgamma} u_N(p)$. Moreover, due to the affine decomposition assumption \eqref{ass: affine decomposition}, the quantities $A_N$ and $f_N$ can be computed via  
\begin{align*}
    A_N(p) = \sum_{q=1}^{Q_a} \theta_a^q(p) V_{\bgamma}^\top A^q V_{\bgamma},\qquad
    f_N(p) = \sum_{q= 1}^{Q_f} \theta_f^{q}(p) V_{\bgamma}^\top f^{q},
\end{align*}
where the matrices $V_{\bgamma}^{\top} A^q V_{\bgamma}$ and vectors $V_{\bgamma}^{\top} f^{q}$ can be pre-computed offline.

Since the linear space $ V_{\bgamma}$ consists of orthonormalized solution vectors for certain parameters, the following question arises: which parameters should we choose to construct the reduced basis? We can formulate this problem as a COP with loss function \eqref{eq: loss-reduced-basis}. Alternatively, we can consider a loss function based on the residual. We define the error function
\begin{align}
\label{eq: residual}
    \err(p, \bgamma) \coloneq \norm{f - A(p) u_{\bgamma}^p}^2_{(A(p))^{-1}} = \norm{f - A(p) \mathbf{V}_{\bgamma} u(p)_N}^2_{(A(p))^{-1}},
\end{align}
where $u_N(p)$ is given by \eqref{eq: reduced system} and $u_{\bgamma}(p) \coloneq V_{\bgamma} u_N(p)$.
In this case, the loss function reads
\begin{align}
\label{eq: loss-residual}
    \calL_{\text{res}}(\bgamma) \coloneq \max_{p \in P} \err(p, \bgamma).
\end{align}

To construct a reduced basis, both GSM and PDM can be applied. In Section~\ref{sec: Numerical examples}, we compare GSM and PDM in a thermal block example.

%% file: Empirical_Interpolation_Method.tex
The Empirical Interpolation Method (EIM) is another technique in Model Order Reduction that utilizes GSM and, hence, also suffers from the scalability issues of GSM. Here, we explain how PDM can be used for EIM. We start by explaining EIM as introduced in \cite{barrault2004}.

In the reduced basis method, one often assumes that the stiffness $A(p)$ and the load vector $f(p)$ admit an affine decomposition. In practice, however, this assumption is not satisfied in general. To overcome this problem, methods such as EIM or hyperreduction are applied to construct an approximate affine decomposition. More precisely, suppose $s(x, p)$ is a non-affine function with $x \in \Omega_h \subset \Omega$ and $p \in P$. The aim is then to express $s(x, p)$ as
\begin{align}
\label{eq: affine approximation}
s(x, p) = \sum_{j=1}^n \beta_j(p) \rho_j(x) + e_{\text{EIM}}(x,p),
\end{align}
where $\sup_{p\in P}\norm{e_{\text{EIM}}(\cdot, p)}_{L^\infty(\Omega_h)} \le \varepsilon$ and $\varepsilon>0$ is a pre-selected tolerance. 

We consider a set of quadrature points $\bX \coloneq \{x_k\,|\, k=1,\ldots, N_q\}$, and we represent the function $s: \Omega_h \times P \to \R$ as a 
vector function $\mathbf{s}: P \to \R^{N_q}$, where 
\begin{align*}
    \mathbf{s}_k(p)\coloneq s(x_k, p), \qquad k = 1, \ldots, N_q.
\end{align*}
Hence, we aim to find a configuration $\bgamma\subset P$ such that $\text{span}\{\mathbf{s}(p)\,|\, p\in\bgamma\}$ approximates the manifold $\{\mathbf{s}(p) \,|\, p \in P\}$.

Given a configuration $\bgamma \in \Gamma_n(P)$, $n \in \N$, we construct an approximation $\mathbf{s}^n(p)$ of $\mathbf{s}(p)$ by
\begin{align}
\label{eq: EIM approximation}
\mathbf{s}^n(p) = \sum_{j=1}^n \beta_j(p) \rho_j.
\end{align}
Here, we need to decide how to determine the vector $\boldsymbol{\beta}(p) = (\beta_1(p), \ldots, \beta_n(p))$. We do this by selecting interpolation points $\bI_n=\{x_\ell\in X\,|\, \ell\in I_n\subset \{1,\ldots,N_q\}\}$ and require that the approximated function and the original function are equal at the interpolation points. We thus require that
\begin{align}
\label{eq: interpolation restriction}
    \mathbf{s}_i^n(p) \coloneq \sum_{j=1}^n \beta_j(p) \rho_j (x_i) = \mathbf{s}_i(p), \quad \ell \in I_n.
\end{align}
To illustrate how $\mathbf{s}^n$ can numerically be computed, we let $\mathbbm{Q}_n \in \R^{N_q \times n}$ be the matrix with columns~$\rho_i$, $i = 1, \ldots, n$. The matrix $\mathbbm{Q}_{I_n} \in \R^{n \times n}$ is then constructed by restricting the matrix $\mathbbm{Q}_n$ to the rows indexed by $I_n$. Here, $I_n$ represents the indices of the selected interpolation points after $n$ steps of the algorithm. Denoting by $\mathbf{s}_{I_n}$ the restriction of $\mathbf{s}$ to the elements of $I_n$, \eqref{eq: interpolation restriction} can then be written as
\begin{align}
    \label{eq: vectorized interpolation restriction}
    \mathbbm{Q}_{I_n} \boldsymbol{\beta}(p) = \mathbf{s}_{I_n}(p).
\end{align}
From this, it follows that $\boldsymbol{\beta}$ can be expressed as 
\begin{align*}
\boldsymbol{\beta}(p) = \big(\mathbbm{Q}_{I_n}\big)^{-1} \mathbf{s}_{I_n}(p).
\end{align*}
Therefore, the approximation $\mathbf{s}^n(p)$ takes the form
\begin{align}
    \label{eq: s_n approximation}
    \mathbf{s}^n(p) = \mathbbm{Q}\mathbbm{Q}_{I_n}^{-1}\mathbf{s}_{I_n}(p).
\end{align}

The question now arises how to determine the functions $\rho_i$ and the interpolation points $I_n$. For the determination of the configuration $\bgamma \in \Gamma_n$, we could apply GSM or PDM with $\err(\mathbf{s}(p), \bgamma) = \norm{\mathbf{s}(p) - \mathbf{s}_n(p)}_{\sup}$. This provides an $\mathbf{s}(p_n)$ at every step $n \in \N$ and we normalize after each step to obtain $\{\rho_1, \ldots, \rho_n\}$. This normalization is based on the interpolation points.

These interpolation points are also determined recursively, where a new interpolation point is selected based on the error estimate in the previous step. In the first step of GSM and PDM, and given the quadrature points $\bX$, we select $\mathbf{s}(p_1)$ such that
\begin{align*}
p_1 = \argmax_{p \in P} \norm{\mathbf{s}(p)}_{\infty},
\end{align*}
if we can determine this value exactly. Otherwise, we select $p_1 \in P$ randomly.
We then set
\begin{align}
\label{eq: first interpolation point}
i_1 \coloneq \argmax_{k = 1, \ldots, N_q} |\mathbf{s}_k(p_1)| \qquad \text{and}\qquad r_1 \coloneq \mathbf{s}_{i_1}(p_1).
\end{align}
Then, we set $\rho_1 \coloneq r_1^{-1}\mathbf{s}(p_1)$.
In the consecutive steps, we select 
\begin{align}
    \label{eq: interpolation point}
    i_{n+1} \coloneq \argmax_{k=1, \ldots, N_q} |\mathbf{s}_k(p_{n+1})-\mathbf{s}_k^n(p_{n+1})| \quad \text{and}\quad r_{n+1} \coloneq \mathbf{s}_{i_{n+1}}(p_{n+1})-\mathbf{s}_{i_{n+1}}^n(p_{n+1}).
\end{align}
In general, $\rho_{n+1}$ is defined as
\begin{align*}
    \label{eq: general rho selection}
    \rho_{n+1} = r_{n+1}^{-1}\big(\mathbf{s}(p_{n+1}) - \mathbf{s}^n(p_{n+1})\big).
\end{align*}
We could prove by induction that, indeed, $\text{span} \{\rho_1, \ldots, \rho_{n}\} = \text{span} \{\mathbf{s}(p_1), \ldots, \mathbf{s}(p_n)\}$.

After convergence of the objective function, $\err(\mathbf{s}(p), \bgamma)$, we have an affine approximation of $\mathbf{s}(p)$ given by $\eqref{eq: EIM approximation}$. In Section \ref{sec: Numerics - subsec: gaussian heat source}, we consider a numerical example where we first apply EIM to find an affine approximation and then use this to find a reduced basis.

%% file: Numerical_examples.tex
In this section, we investigate two numerical examples. The first example is a thermal block model \cite{rozza2008reduced}, where the underlying PDE satisfies the affine decomposition assumption. We compare the performance of PDM with GSM-Random and GSM-LHS.

The second example is a model for heat conduction with a Gaussian heat source, where the affine decomposition is not satisfied. Therefore, we first apply EIM both with GSM-Random and PDM on the right-hand side of the PDE and then use the approximated affine decomposition to construct a reduced basis.\footnote{The code for PDM used in this manuscript can be found here:
https://doi.org/10.5281/zenodo.13935184
The collection of all the codes used to generate the numerical results presented in Section~\ref{sec: Numerical examples} can be found here:
https://gitlab.tue.nl/s158446/polytope-division-paper}

\subsection{Thermal block example}
\begin{figure}
    \centering
    \begin{minipage}[b]{0.49\textwidth}
  \centering
  \input{Figures/ThermalBlock-6d}
  \vspace{2.5em}
  \centering
  \caption{Thermal Block with $6$ parameters}
  \label{fig: tb6} 
    \end{minipage}
    \begin{minipage}[b]{0.49\textwidth}
  \includegraphics[width=\linewidth]{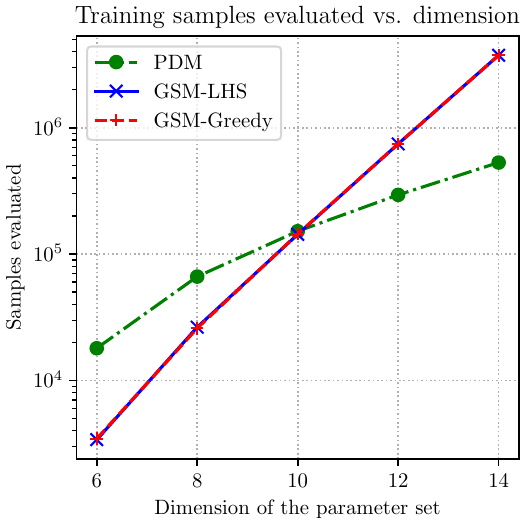}
  \caption{Samples evaluated vs.\ dimension}
  \label{fig: TB samples}
    \end{minipage}
\end{figure}

In the thermal block example, we consider the following stationary heat equation
\begin{align*}
    \begin{cases}
    - \nabla \cdot [\kappa(x, p) \nabla u(x, p)] = 1 \qquad &\textup{for } x \in \Omega,\\
     u(x, p) = 0 \qquad &\textup{for } x \in \partial \Omega.
    \end{cases}
\end{align*}
Here, $\Omega = [0,1]^2$ is divided into blocks, $B_i$, $i = 1, \ldots, d$, where $d$ is the number of parameters. Each block has thermal conductivity  $p_i$ for $i \in \{1, \ldots, d\}$ and $\kappa(p, x) = p_i$ if $x \in B_i$. For all $i \in \{1, \ldots, d\}$, we have that $1 \le p_i \le 10$. We consider thermal blocks with different numbers of parameters with each block having $2$ rows and $m$ columns with $m \in \{3, \ldots, 7\}$. In this way, we can compare the computation times of GSM-Random, GSM-LHS, and PDM for different dimensions of the parameter space. Figure~\ref{fig: tb6} depicts a thermal block with $6$ parameters.

We apply GSM-Random, GSM-LHS, and PDM to find a reduced basis for the solution of the thermal block problem in the finite element space $H_0^1(\Omega)$ with $2371$ degrees of freedom. The parameter space is given by the hypercube, $\calP = [1, 10]^d$. For GSM-Random and GSM-LHS, we choose the sample size $\samsize = 2^d$, where $d$ is the dimension. The reduced basis is computed with tolerance $\eps = 10^{-12}$, and we compare the total number of evaluated training samples to construct this basis vs. the dimension in Figure \ref{fig: TB samples}. We note that fewer samples are evaluated in GSM-Random and GSM-LHS for a lower dimension of the parameter set. However, in higher dimensions, the number of training samples evaluated for PDM is much lower than for GSM-Random and GSM-LHS, leading to lower computation times.

We should remark here that this comparison depends strongly on the sample size $\samsize$ that we choose for GSM-Random and GSM-LHS. However, in practice, the user also has to choose the sample size $\samsize$ without knowing a priori what sample size is sufficient. Choosing $2^d$ samples means roughly taking $2$ samples per dimension, which could be considered very low---clearly, one does not take only $4$ sample points for $d=2$. In this sense, PDM substantially differs from both versions of GSM because the user does not have to choose a sample size a priori. Instead, the algorithm keeps refining the sample space until the error estimates in all the barycenters are below the desired tolerance.

\begin{figure}[h]
  \centering
  \includegraphics[width=.44\linewidth]{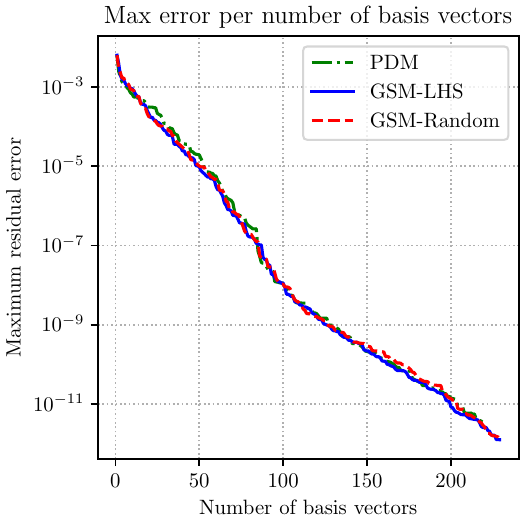}
  \caption{Verification stage: Maximum squared residual error for the reduced bases of the thermal block problem for PDM, GSM-LHS, and GSM-Random}
  \label{fig: TB verification}
\end{figure}

We also compare the quality of the reduced bases constructed by GSM-Random, GSM-LHS, and PDM for dimension $d = 14$. We do this by considering a set of verification samples and computing the squared residual error \eqref{eq: residual}. We then select the maximum error of all the verification samples per number of basis vectors. Figure \ref{fig: TB verification} displays this result. In this case, we took $5000$ verification samples. We note here that the squared maximum residual error is of the same order of magnitude for GSM-Random, GSM-LHS, and PDM, despite the large difference in number of samples evaluated shown in Figure~\ref{fig: TB samples}.

\begin{figure}[h]
    \centering
    \begin{minipage}[b]{0.49\linewidth}
        \centering
  \includegraphics[width=\linewidth]{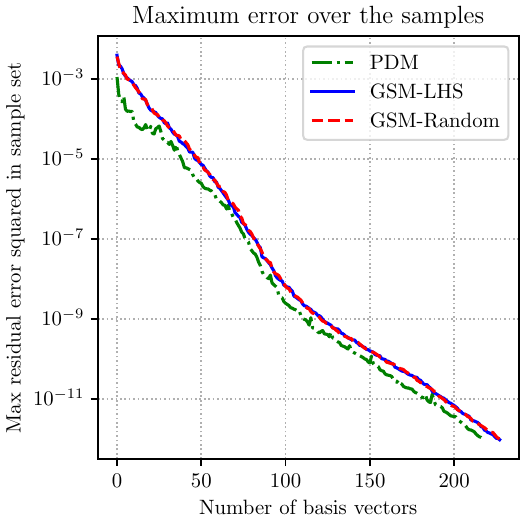}
  \caption{Maximum squared residual error over the training samples for the reduced basis approximation for the thermal block problem.}
  \label{fig: TB convergence plots combined}
    \end{minipage}
    \hfill
    \begin{minipage}[b]{0.49\linewidth}
  \centering
  \includegraphics[width=\linewidth]{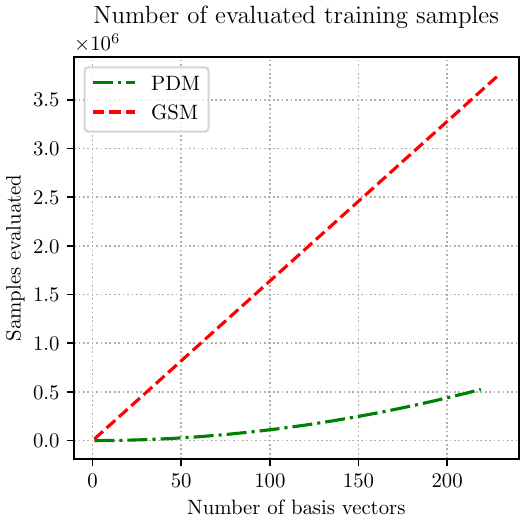}
  \caption{Comparison of the evaluated number of training samples in the $14$-dimensional thermal block per number of basis vectors.}
  \label{fig: TB evaluated samples}
\end{minipage}
\end{figure}

Moreover, in Figure \ref{fig: TB convergence plots combined}, we investigate the convergence of the maximum error in all the training sample points for GSM-Random and GSM-LHS, and all the barycenters in PDM for the case $d=14$. In Figure \ref{fig: TB convergence plots combined}, we observe that the convergence rates are comparable. One difference is that the squared maximum residual error in the training sample set is monotonically decreasing for GSM-Random and GSM-LHS. In contrast, the error in the barycenter set of PDM oscillates slightly around $25$ basis vectors. We postpone the discussion of this behaviour to the next example.

Lastly, Figure \ref{fig: TB evaluated samples} displays the number of samples evaluated per number of basis vectors. We disregarded the graph for GSM-LHS due to its similarity with that of GSM-Random. The figure clearly depicts a large difference in the number of evaluated samples of both GSM methods and PDM. This difference also leads to a major difference in computational costs.

\subsection{Heat conduction model with Gaussian heat source}
\label{sec: Numerics - subsec: gaussian heat source}
We recall from Section~\ref{subsec: Empirical Interpolation Method} that the Empirical Interpolation Method (EIM) is used when the affine decomposition assumption fails, and it is another example of COP.
In this section, we explore an example where the affine decomposition assumption is not satisfied, applying EIM with both GSM-Random and PDM. We apply GSM-Random and PDM twice in this example. First, we apply both methods to EIM to get an approximated affine decomposition with tolerance $\eps = 10^{-6}$. Second, we apply the methods to compute a reduced basis with tolerance $\eps = 10^{-10}$.

For this example, we consider a modification of EIM for a non-affine elliptic problem based on the Tutorial from \cite{RBnicsCite}. In particular, we consider a heat conductivity model with a Gaussian heat source. The conductivity coefficient in this model always equals $1$ and $u(x, p)$ is the temperature at position $x \in \Omega = (-1,1)^2$ for parameter 
\[
    p = (\mu_1, \mu_2, \sigma_1, \sigma_2, \rho) \in P\coloneq [-1,1]^2 \times [1,3]^2 \times [-0.8, 0.8] \subset \R^5.
\]
The function $u$ is then a solution of the following PDE
\begin{align*}
    -\Delta u = g(\cdot, p) \quad \text{in $\Omega$},\qquad u(\cdot,p) = 0\quad\text{on $\partial\Omega$},
\end{align*}
where the heat source $g$ is given by
\begin{align*}
    g(x, p) = N(p) \cdot \exp\left\{-\frac{1}{2(1-\rho^2)}\Bigg[\Big(\frac{x_1-\mu_1}{  \sigma_1}\Big)^2 - 2\rho \frac{x_1-\mu_1}{\sigma_1}\frac{x_2 - \mu_2}{\sigma_2} + \Big(\frac{x_2-\mu_2}{\sigma_2}\Big)^2\Bigg]\right\},
\end{align*}
where $N(p)$ is a normalization factor given by
\begin{align*}
    N(p) = \frac{1}{2\pi \sigma_1 \sigma_2 \sqrt{1-\rho^2}}.
\end{align*}

As indicated in the thermal block example, the computational costs of the offline stage of the GSM-Random approach depend strongly on the chosen sample size $\samsize$. In the previous example, we considered a random sample set of size $\samsize = 2^d$, roughly corresponding to $2$ samples per dimension, as is often done in practice for large dimensions. 

\begin{figure}[h]
\centering
\begin{subfigure}[b]{0.49\linewidth}
\centering
{\includegraphics[width=\linewidth]{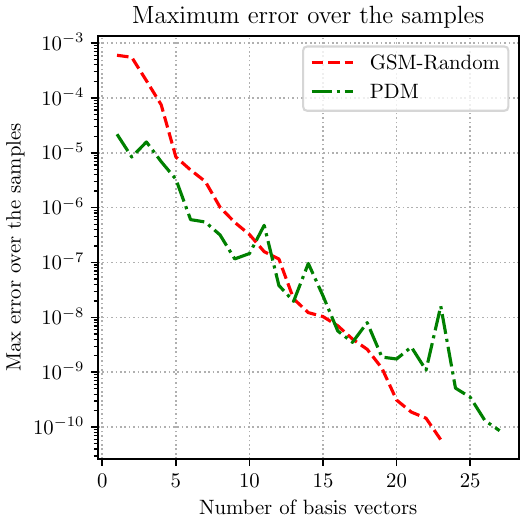}}
\caption{}
\label{fig: sf: convergence plots combined}
\end{subfigure}
\hfill
\begin{subfigure}[b]{0.49\linewidth}
\centering
{\includegraphics[width=\linewidth]{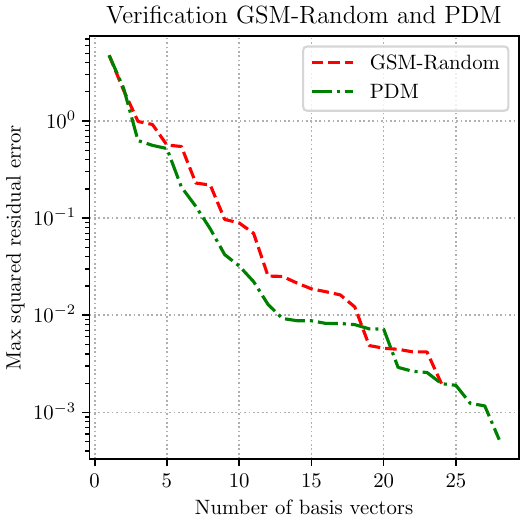}}
\caption{}
\label{fig: sf: verification EIM}
\end{subfigure}
\caption{(a) Maximum error over the $\samsize=481$ training samples per number of basis vectors. (b) Verification of the reduced basis for a set of $1000$ verification samples.}
\label{fig: EIM}
\end{figure}
To get a richer comparison, we choose the sample size differently in this example. Since we do not have to determine a sample size a priori for PDM, we first perform this algorithm to get an affine approximation via EIM. We then use the result of PDM-EIM to compute a reduced basis and consider the final number of barycenters.
For PDM-EIM, the final number of samples is $\samsize=481$ resulting in $183$ elements in the reduced basis.

Next, we perform GSM-EIM with a sample size of $481$. We then use the result of GSM-EIM to perform the reduced basis method with $183$ samples. We note that GSM-Random still evaluates more samples than PDM overall , since the number of samples in PDM increases at each step, while it remains constant in GSM-Random.

We then compare the numerical performance of the two methods by plotting the maximum squared residual error in all the GSM-Random sample points and all the PDM barycenters, as shown in Figure \ref{fig: sf: convergence plots combined}.
For GSM-Random, the maximum squared residual error decreases monotonically, while for PDM it oscillates. This is because GSM-Random uses a fixed sample set, leading (in this case) to a monotone decrease in the objective function after each update. In contrast, PDM refines the polytope division at each step, often targeting areas with the highest expected average error, which can reveal new regions where the objective function is larger; therefore, we cannot expect a monotone decrease even in this simple case.

We also compare the numerical performance of the reduced bases constructed with GSM-Random and PDM in a set of verification samples. We consider $1000$ verification samples and compute the maximum squared residual error versus the number of basis vectors. Figure \ref{fig: sf: verification EIM} displays this result. The reduced basis constructed with GSM-Random has fewer basis vectors than PDM because the algorithm stops after the tolerance is attained and this happens at an earlier stage for GSM-Random.

\begin{figure}[h]
\centering
\begin{minipage}[b]{0.44\linewidth}
{\includegraphics[width=\linewidth]{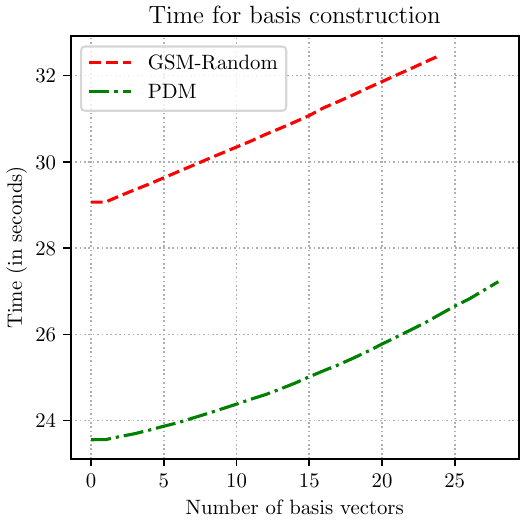}}
\subcaption{}
\label{fig: sf: eim-time}
\end{minipage}
\hfill
\begin{minipage}[b]{0.44\linewidth}
{\includegraphics[width=\linewidth]{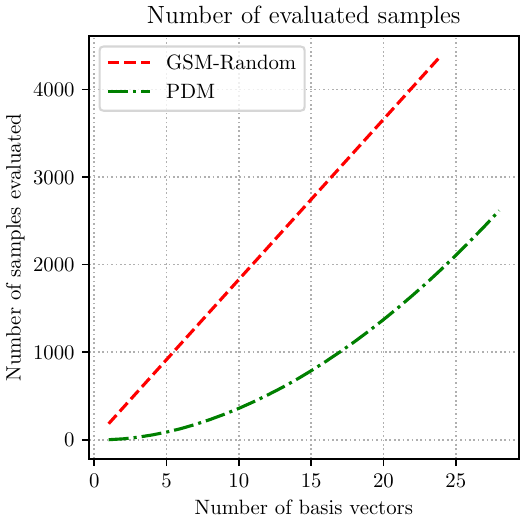}}
\subcaption{}
\label{fig: sf: evaluated samples}
\end{minipage}
\caption{Comparison of computation times and the evaluated number of samples of GSM-Random and PDM in the Gaussian heat source example per number of basis vectors}
\label{fig: eim-time-plots}
\end{figure}

In Figure \ref{fig: eim-time-plots}, we compare the computation times of GSM-Random and PDM and evaluate the total number of error estimates computed. In Figure \ref{fig: sf: eim-time}, the intersections of the graph with the $y$-axis represent the time to perform GSM-EIM and PDM-EIM, respectively on a computer with 2.6 GHz Intel Core i7 processor and 32 GB of RAM. We note that PDM is faster than GSM-Random, even though the errors in the verification sets are lower for PDM. The idea behind PDM is that we adaptively refine to find the relevant samples. This experiment seems to indicate that PDM finds more relevant samples than GSM-Random in considerably less time. In Figure \ref{fig: sf: evaluated samples}, we display the total number of samples for which the error estimate is computed vs the number of basis vectors. This difference in number of evaluated samples explains the difference in computation times.

%% file: Figures/ThermalBlock-6d.tex
\begin{tikzpicture}[scale=1.0]
        \def\d{5};
        \def\o{2};

        \coordinate (P11) at (0,0);
        \coordinate (P12) at (\d/3, 0);
        \coordinate (P13) at (2*\d/3, 0);
        \coordinate (P14) at (\d,0);

        \coordinate (P21) at (0,\d/2);
        \coordinate (P22) at (\d/3, \d/2);
        \coordinate (P23) at (2*\d/3, \d/2);
        \coordinate (P24) at (\d,\d/2);

        \coordinate (P31) at (0,\d);
        \coordinate (P32) at (\d/3, \d);
        \coordinate (P33) at (2*\d/3, \d);
        \coordinate (P34) at (\d,\d);

        \coordinate (B1) at (\d/6, 3*\d/4);
        \coordinate (B2) at (\d/2, 3*\d/4);
        \coordinate (B3) at (5*\d/6, 3*\d/4);

        \coordinate (B4) at (\d/6, \d/4);
        \coordinate (B5) at (\d/2, \d/4);
        \coordinate (B6) at (5*\d/6, \d/4);

        \draw[draw=none, fill=red, opacity=.5] (P11) -- (P12) -- (P22) -- (P21);
        \draw[draw=none, fill=blue, opacity=.35] (P12) -- (P13) -- (P23) -- (P22);
        \draw[draw=none, fill=yellow, opacity=.5] (P13) -- (P14) -- (P24) -- (P23);
        \draw[draw=none, fill=green, opacity=.5] (P21) -- (P22) -- (P32) -- (P31);
        \draw[draw=none, fill=orange, opacity=.5] (P22) -- (P23) -- (P33) -- (P32);
        \draw[draw=none, fill=cyan, opacity=.5] (P23) -- (P24) -- (P34) -- (P33);

        \draw[thick, color=black] (P12) -- (P32);
        \draw[thick, color=black] (P13) -- (P33);
        \draw[thick, color=black] (P21) -- (P24);
        \draw[thick, color=black] (P11) -- (P14) -- (P34) -- (P31) -- (P11);

        \node at (B1) {$B_1$};
        \node at (B2) {$B_2$};
        \node at (B3) {$B_3$};
        \node at (B4) {$B_4$};
        \node at (B5) {$B_5$};
        \node at (B6) {$B_6$};
\end{tikzpicture}

%% file: Outlook.tex
In this paper, we introduced the Polytope Division Method. PDM is an algorithm that can be used to find a heuristic solution to Configuration Optimization Problems. The strength of PDM over the Greedy Sampling Method is twofold:
\begin{enumerate}[label={(\arabic*)}]
    \item It requires significantly fewer sample evaluations in high dimensions;
    \item It does not require a pre-defined training sample size.
\end{enumerate}

PDM could be further explored for various practical applications and theoretical convergence results. In future work, we aim to provide a rigorous foundation for PDM and consider further applications of PDM to COPs, e.g., in the context of context of Optimal Experimental Design (OED) and Active Learning for Regression (ALR).

%% file: Proofs.tex
Here, we provide the proof of Lemma~\ref{lemma: facets p-boundary polytopes}. The proof relies on the following intermediate result:
\begin{lemma}
\label{lemma: conv_x_F1 facets}
    Let $P \subset \R^d$ be a polytope. Let $F \in \partial P$ be a facet of $P$ and $x \in \mathrm{int}(P)$ be an interior point of $P$. Consider the polytope $\tau \coloneq \conv(x \cup F)$. Then,
    \[
        \partial\tau = \{F\}\cup \bigl\{\conv(x \cup G):  G \in \partial F\bigr\}.
    \]
\end{lemma}
\begin{proof}
Since $x\in \rm{int}{P}$, it is not difficult to see that $F$ is indeed a facet of $\tau$. We are then left to prove that $\conv(x \cup G)$ is also a facet of $\tau$.
A facet of $\tau$ needs to be a $(d-1)$-dimensional polytope. We note that $F$ is a $(d-1)$-dimensional polytope and $G$ is $(d-2)$-dimensional.
Since $x \in \rm{int}(P)$, it is not coplanar to the vertices of $F$ and is therefore not coplanar to the vertices of $G$. Hence, $\conv(x \cup G)$ is a $(d-1)$-dimensional polytope. In particular, there exist some  $c_1 \in \R^d$ and $\alpha \in \R$ such that
\begin{align*}
\conv(x \cup G) \cap \calH = \conv(x \cup G),\qquad\text{where $\calH \coloneq \{y \in \R^d\, |\, c_1 \cdot y = \alpha \}$},
\end{align*}
i.e., $\calH$ is a hyperplane through $G$. We show next that $\calH$ is a separating hyperplane.
First, we show $\tau \cap \calH = \conv(x \cup G)$. The inclusion $\tau \cap \calH \supset \conv(x \cup G)$ is trivial.
Now let $y \in \tau \cap \calH$. We note that $y \in \tau$ implies there exist $\lambda_1, \ldots, \lambda_{k+1} \ge 0$, $\sum_{i=1}^{k+1} \lambda_i = 1$ such that
$$y = \sum_{i=1}^k \lambda_i v_i + \lambda_{k+1} x.$$
Here, $v_1, \ldots, v_k$ are the vertices of $F$. Without loss of generality, we say $v_1, \ldots, v_l$, are the vertices of $G$ for some $l < k$. If for all $i \in \{l+1, \ldots, k\}$ it holds that $\lambda_i = 0$ then $y \in \conv(x \cup G)$. But if there exists at least one $i \in \{l+1, \ldots, k\}$ such that $\lambda_i > 0$, then $\calH$ must span a $d$-dimensional space, which contradicts the fact that $\calH$ is a hyperplane. Hence, we conclude that $y \in \conv(x \cup G)$.

Next, we show that for all $v \in \{v_{l+1},\ldots, v_k\}$, it holds that $c_1 \cdot v > \alpha$ (or $c_1 \cdot v < \alpha)$. In other words, we show that the vertices of $F$ not contained in $G$ lie on one side of the hyperplane $\calH$. Since $G \in \partial F$, there exist $c_2 \in \R^d$ and $\beta \in \R$ such that $G$ is contained in the hyperplane $\calH_G = \text{span}\{ y\in\R^d: c_2\cdot y=\beta\}$ and $c_2 \cdot y > \beta$ for all $y \in \{v_{l+1},\ldots, v_k\}$. For the sake of contradiction, suppose there exist $v, w \in \{v_{l+1},\ldots, v_k\}$ such that $c_1 \cdot v > \alpha$ and $c_1 \cdot w < \alpha$.
If such $v, w$ exist, there is a $\lambda \in (0,1)$ such that $c_1 \cdot (\lambda v + (1-\lambda) w) = \alpha$.
Hence, $\lambda v + (1-\lambda) w \in F \subset \tau$ and $\lambda v + (1-\lambda) w \in \calH$, implying that $\lambda v + (1-\lambda) w \in G$, i.e.,
$$
    c_2 \cdot (\lambda v + (1-\lambda) w) = \beta.
$$
On the other hand, since $v,w$ satisfies $c_2\cdot v >\beta$ and $c_2\cdot w >\beta$, we deduce
\begin{align*}
c_2 \cdot (\lambda v + (1-\lambda) w) &= \lambda c_2 \cdot v + (1-\lambda) c_2 \cdot w > \lambda \beta + (1-\lambda) \beta = \beta,
\end{align*}
leading to a contradiction. Hence, without loss of generality, we suppose that $c_1 \cdot v > \alpha$ for all $v \in \{v_{l+1},\ldots, v_k\}$.
But then $\calH$ is a separating hyperplane of $\tau$. Indeed, for any $y \in \tau$, it holds that there exist $\mu_1, \ldots, \mu_{k+1}$ such that $\sum_{i=1}^{k+1} \mu_i =1$ and $y = \sum_{i=1}^k \mu_i v_i + \mu_{k+1} x$. Then,
\begin{align*}
c_1 \cdot y = \sum_{i=1}^k \mu_i c_1 \cdot v_i + \mu_{k+1} c_1 \cdot x \ge \sum_{i=1}^{k+1} \mu_i \alpha = \alpha.
\end{align*}
This proves that $\conv(x \cup G)$ is indeed a facet of $\tau$.

Lastly, we have to prove that $\tau$ does not have any other facets. We note that $\tau$ is the convex hull of $\{x, v_1, \ldots, v_k\}$ with $v_1, \ldots, v_k$ the vertices of $F$. All possible facets can be described by the convex hull of a subset of these vertices that describe a $(d-1)$-dimensional polytope.
Since $F$ is a facet, there exists a separating hyperplane $\calH$ such that
$\tau \cap \calH = F$.
Because $F$ is $(d-1)$-dimensional, this hyperplane is unique, meaning that if there exists a subset of vertices $\{w_1, \ldots, w_k\} \subset \partial_{d-1} F$ such that $\calN=\text{span}(\{w_1, \ldots, w_k\})$ and $\calN$ is $(d-1)$-dimensional, then $\calN = \calH$. Hence, the convex hull of a strict subset of vertices of $F$ does not describe a facet.
The only other candidate for a facet is $\conv(\{x, w_1, \ldots, w_k\})$, where $w_i \in \partial_{d-1} F$ for $i \in \{1, \ldots, k\}$, such that $\conv(\{w_1, \ldots, w_k\})$ is $(d-2)$-dimensional.
But then $\conv(\{w_1, \ldots, w_k\})$ has to be a facet of $F$, for otherwise, there would be no separating hyperplane $\calH$ such that
$$F \cap \calH = \conv(\{w_1, \ldots, w_k\}).$$
In particular, the hyperplane through the points $x, w_1, \ldots, w_k$, is not a separating hyperplane. Hence, we conclude that all the facets are given by $F$ and $\conv(x \cup G)$ with $G \in \partial F$.
\end{proof}

We conclude the appendix with the proof of Lemma~\ref{lemma: facets p-boundary polytopes}.

\begin{proof}[Proof of Lemma~\ref{lemma: facets p-boundary polytopes}]
We prove this by induction over the dimension $d\ge 2$.
    For $d = 2$, we only have to consider the object $\conv(x_1 \cup F)$, where $F \in \partial P$. The statement then follows directly from Lemma~\ref{lemma: conv_x_F1 facets}.
    
    Our induction hypothesis states that for a $d$-dimensional hypercube $P$ and a $P$-boundary polytope of the form $\tau = \conv(\{x_1, \ldots, x_n\} \cup F)$, with $F \in \partial_n P$, it indeed holds that
    \begin{align*}
        \partial \tau = \bigcup_{i = 1}^n \conv\left(\{x_1, \ldots, x_n\} \backslash \{x_i\} \cup F \right) \cup \bigcup_{G \in \partial F} \conv(\{x_1, \ldots, x_n\} \cup G).
    \end{align*}
    Now, we prove it is also true for dimension $d+1$. 
    
    Let $Q$ be a $(d+1)$-dimensional hypercube and $E\in\partial_n Q$. We first consider the case for $n = 1$ separately, i.e., consider the $(d+1)$-dimensional polytope $\conv(x_1 \cup E)$, then according to Lemma \ref{lemma: conv_x_F1 facets}, the facets are given by $E$ and $\conv(x_1 \cup G)$, where $G \in \partial E$. Hence, the statement is true for dimension $d+1$ and the case $n=1$.
    
    Consider a $Q$-boundary polytope $\varrho = \conv(\{x_1, \ldots, x_n\}\cup E)$, with $n \in \N$ and $2 \le n \le d$ and $E \in \partial_n Q$. Since $\varrho$ is a $(d+1)$-dimensional polytope, there exists a polytope $R$ for which $L\coloneq\conv(\{x_1, \ldots, x_{n-1}\} \cup E)\in \partial R$ and $x_n\in \interior{R}$. Hence, by Lemma \ref{lemma: conv_x_F1 facets},
    \[
        \partial \varrho = \{L\} \cup \bigl\{ \conv(x_n \cup G) : G\in \partial L\bigr\}.
    \]
    Therefore, to determine all the facets of $\varrho$, we have to determine all the facets of $L$. This latter object is $d$-dimensional and $1 \le n-1 \le d$. Hence, by exploiting the induction hypothesis, we have that 
    \begin{align*}
        \partial L = \bigcup_{i = 1}^{n-1} \conv\bigl(\{x_1, \ldots, x_{n-1}\} \backslash \{x_i\} \cup E \bigr) \cup \bigcup_{K \in \partial E} \conv(\{x_1, \ldots, x_{n-1}\} \cup K)
    \end{align*}
    Consequently,
    \begin{align*}
        \bigcup_{G\in \partial L} \conv(x_n \cup G) = \bigcup_{i = 1}^{n-1} \conv\left(\{x_1, \ldots, x_n\} \backslash \{x_i\} \cup E \right) \cup \bigcup_{K \in \partial E} \conv(\{x_1, \ldots, x_n\} \cup K).
    \end{align*}
    Therefore, we can conclude
    \begin{align*}
        \partial \varrho = \bigcup_{i = 1}^{n} \conv\left(\{x_1, \ldots, x_n\} \backslash \{x_i\} \cup E \right) \cup \bigcup_{G \in \partial E} \conv(\{x_1, \ldots, x_n\} \cup G).
    \end{align*}
    Hence, the statement is also true for dimension $d+1$. 
\end{proof}